\def\WHO{nbd}
\def\version{10.7.2021}\def\users{}  %
\RequirePackage{fix-cm}
\documentclass[global,final,numbook,envcountsame]{article}                     % onecolumn (standard format)
\usepackage{graphicx}

\usepackage{amsmath}
\numberwithin{equation}{section}

\usepackage{authblk}
\usepackage{upgreek}

\usepackage{fullpage}
\usepackage{bm}
\usepackage{amsmath}
\usepackage{hyperref}
\usepackage{amsfonts}
\usepackage{amssymb}
\usepackage{xcolor}

\usepackage{mathrsfs} % THIS PACKAGE ALLOWS FOR ``scr'' FONTS

\def\users{final-layout}  
\usepackage{ifthen}
\ifthenelse{\equal{\users}{final-layout}}{}{
\usepackage{fancyhdr}
\pagestyle{fancy}
\headheight=28pt\headwidth=17cm
\definecolor{gray}{gray}{0.5}
\rhead{\color{gray}geophysical models\\
T.Roub\'\i\v cek \& G. Tomassetti}
\chead{}
\lhead{Version \version, file: \jobname.tex, \underline{\Large\color{red}\WHO's working on} \\
compiled:
\number\day.\number\month.\number\year\ at
\the\hour:\ifnum\minute<10 0\fi\the\minute\ h }
}

\usepackage[notcite,notref,color]{showkeys}
\definecolor{labelkey}{rgb}{1.,.2,0.}

\newcount\hour \newcount\minute
\hour=\time
\divide \hour by 60
\minute=\time
\loop \ifnum \minute > 59 \advance \minute by -60 \repeat

\usepackage[normalem]{ulem}

\usepackage[normalem]{ulem}
\usepackage{ifthen}
\usepackage{color}

\ifthenelse{\equal{\users}{final-layout}}{
	
	\newcommand{\INSERT}[1]{#1}
	\newcommand{\COMMENT}[1]{}
	\newcommand{\COMMENTGT}[1]{}
	\newcommand{\TODO}[1]{}
	\newcommand{\INTERNAL}[1]{}
	\newcommand{\QUESTION}[1]{}
	\newcommand{\DELETE}[1]{}

	\newcommand{\REM}[1]{\marginpar{\bfseries\tiny{\color{blue}}}}
    \newcommand{\MARGINOTE}[1]{}
}
{
	\newcommand{\INSERT}[1]{{\color{blue}\uuline{#1}\color{black}}}
	\newcommand{\COMMENT}[1]{{\color{red}\uuline{#1}\color{black}}}
	\newcommand{\COMMENTGT}[1]{{\hfill\large\color{red}***{#1}***\color{black}\hfill}\\}
	\newcommand{\TODO}[1]{{\color{red}\uuline{#1}\color{black}}}
	\newcommand{\INTERNAL}[1]{\footnote{#1}}
	\newcommand{\QUESTION}[1]{{\color{brown}\uuline{#1}\color{black}}}
	\newcommand{\DELETE}[1]{{\color{red}\sout{#1}\color{black}}}

	\newcommand{\REM}[1]{\marginpar{\bfseries\tiny{\color{blue}#1}}}
\newcommand{\MARGINOTE}[1]{\marginpar{\color{red}\tiny\texttt{#1}}}
}
\def\vv{{\bm v}}

\def\zz{\chi}
  \def\vvk{\vv_\tau^k}
  \def\vvkk{\vv_\tau^{k-1}}
  \def\zzk{\zz_\tau^k}
  \def\zzkk{\zz_\tau^{k-1}}
\def\alphak{\alpha_\tau^k}
  \def\alphakk{\alpha_\tau^{k-1}}
  \def\Eek{{\bm E}_{\tau}^k}
  \def\Eekk{{\bm E}_{\tau}^{k-1}}
 \def\Eetau{{\bm E}_{\tau}^{}}
\def\overlineEetau{\hspace*{.2em}\overline{\hspace*{-.2em}\bm E}_{\tau}^{}}
\def\overlineStau{\hspace*{.2em}\overline{\hspace*{-.2em}\bm S}_{\tau}^{}}
\def\Eptau{{\Ep}_{\tau}^{}}
\def\overlineEptau{\hspace*{.2em}\overline{\hspace*{-.2em}\Ep}_{\tau}^{}}

 \def\Epk{{\Ep}_{\tau}^k}
 \def\Epkk{{\Ep}_{\tau}^{k-1}}
 \def\muk{\mu_\tau^k}
\newcommand\density{\varrho}

\newcommand\DELETEDELETE[1]{}
\usepackage{amsthm}
\theoremstyle{remark}
\newtheorem{theorem}{Theorem}
\newtheorem{remark}{Remark}
\newtheorem{definition}{Definition}
\newtheorem{lemma}{Lemma}

\newcommand\DT[1]{\mathchoice
                 {{\buildrel{\hspace*{.1em}\text{\LARGE.}}\over{#1}}}
                 {{\buildrel{\hspace*{.1em}\text{\Large.}}\over{#1}}}
                 {{\buildrel{\hspace*{.1em}\text{\large.}}\over{#1}}}
                 {{\buildrel{\hspace*{.1em}\text{\large.}}\over{#1}}}}
\newcommand\mtd[1]{\frac{{\rm D}#1}{{\rm D}t}}
\newcommand\mtdk[1]{\frac{{\rm D}_{k-1}^k#1}{{\rm D}_\tau t}}

\newcommand\pdt[1]{\frac{\partial{#1}}{\partial t}} %Partial Derivative with respect to T
\newcommand\Ee{{\bm E}}              %Eulerian elastic strain
\newcommand\Ep{{\bm\varPi}}              %Eulerian plastic strain
\newcommand\mtdEe{\mtd{{\bm E}}}           
\newcommand\mtdEp{\mtd{{\Ep}}}

\newcommand\mtdkEe{\mtdk\Ee}           
\newcommand\mtdkEp{\mtdk\Ep}

\def\widetildeEptau{\hspace*{.1em}\widetilde{\hspace*{-.1em}\Ep}_{\tau}^{}}
\def\widetildeEp{\hspace*{.1em}\widetilde{\hspace*{-.1em}\Ep}}
\newcommand\R{\mathbb{R}}

\renewcommand\d{{\rm d}}
\newcommand\EE{{\bm e}}
\newcommand{\lineunder}[2]{\LU{\begin{array}[t]{c}\underbrace{#1}\vspace*{.5em}\end{array}}{\mbox{\footnotesize\rm #2}}}
\newcommand{\LU}[2]{\begin{array}[t]{c}#1\vspace*{-1em}\\_{#2}\end{array}}
\newcommand{\linesunder}[3]{\LSU{\begin{array}[t]{c}\underbrace{#1}\vspace*{.5em}\end{array}}{\mbox{\footnotesize\rm #2}}{\mbox{\footnotesize\rm #3}}}
\newcommand{\LSU}[3]{\begin{array}[t]{c}#1\vspace*{-1em}\\_{#2}\vspace*{-.5em}\\_{#3}\end{array}}

\newcommand{\Vdots}{\mathchoice{\,\vdots\,}{\:\begin{minipage}[c]{.1em}\vspace*{-.4em}$^{\vdots}$\end{minipage}\;}
{\:\tiny\vdots\:}{\:\tiny\vdots\:}}
\def\FFF{\color{black}}
\def\SSS{\color{black}}
\def\TTT{\color{black}}

\def\EEE{\color{black}}

\title{A convective model for poro-elastodynamics with damage and fluid flow towards
Earth lithosphere
modelling.}

\author[1,2]{Tom\'a\v{s} Roub\'\i\v cek}
\author[3]{Giuseppe Tomassetti}

\affil[1]{ Mathematical Institute, Charles University, \protect\\
  Sokolovsk\'a~83, CZ-186~75~Praha~8, Czech Republic.\protect\\\mbox{}}
  \affil[2]{
 Institute of Thermomechanics, Czech Academy of Sciences,\protect\\
    Dolej\v skova 5, CZ-182 00 Praha 8, Czech Republic.\protect\\
              Email: \texttt{tomas.roubicek@mff.cuni.cz} \protect\\\mbox{}}

\affil[3]{ Department of Engineering, Roma Tre University,\protect\\ Via Vito Volterra 62,  00146 Roma, Italy.\protect\\
           Email: \texttt{giuseppe.tomassetti@uniroma3.it}}
\date{April 17, 2021}
\begin{document}

\maketitle

\begin{abstract}
 Devised towards geophysical applications for various processes in \SSS the
 \EEE lithosphere \TTT or the crust\EEE, a model of poro-elastodynamics with
inelastic
strains and other internal variables like damage (aging) and porosity as
well as with diffusion of water is formulated fully in the Eulerian setting.
Concepts of gradient of the total strain rate as well as the additive
splitting of the total strain rate are used while eliminating the displacement
from the formulation. It relies on that the
elastic strain is small while only the inelastic and the total strains can be
large. The energetics behind this model is derived and used for analysis as far
as the existence of global weak energy-conserving solutions concerns. By this
way, the model in [V. Lyakhovsky et al., Pure Appl. Geophys., 171:3099--3123, 2014] and 
[V. Lyakhovsky et al., Izvestiya, Physics of the Solid Earth, 43:13--23, 2007] is completed to make it mechanically consistent and amenable for analysis.
\end{abstract}

\noindent\textbf{AMS subclass.}{
 35K87, % Systems of parabolic variational inequalities
 35Q74, % PDEs in connection with mechanics of deformable solids
% 35Q79, % PDEs in connection with classical thermodynamics and heat transfer
35Q86, %PDEs in connection with geophysics
%74A15, % Thermodynamics
74A30, % Nonsimple materials
%74A45, % Theories of fracture and damage%
74C10, % Plastic materials - internal variable - small-strain, rate-dependent theories
74C20, % Plastic materials - internal variable - Large-strain, rate-dependent theories
74H20, % Dynamical problems, Existence of solutions
74L05, % Geophysical solid mechanics
74R20, % Anelastic fracture and damage
%80A20, % Heat and mass transfer, heat flow
76S05, % fluid mechanics, Flows in porous media
86A17. %Global dynamics, earthquake problems}
}

\section{Introduction}
%        ~~~~~~~~~~~~
Geophysical models of \SSS the \EEE solid \TTT Earth (i.e.\ particularly the
lithosphere and the crust) \EEE are extremely challenging
applications of continuum mechanics. Such models should capture a lot of
phenomena on various time-space scales
and usually are focused on only specific aspects, cf.\
\cite[Fig.\,1]{BeZi08CBEF} for the spatiotemporal scales relevant for
earthquakes and fault dynamics. On short time scales, fast rupture of
\SSS lithospheric \EEE faults, tectonic earthquakes,
and seismic waves are most prominent phenomena. On large time scales,
aseismic creep, healing of damaged faults, and water (or sometimes oil)
transport in porous rocks are dominant effects to capture.
The water transport processes are intimately coupled with mechanical properties
and possibly also with evolution of porosity and of damage (called also aging
in geophysical applications). Other effects would be heat production and
transfer, magnetism, or volcanism, but
% , except Remark~\ref{rem-thermo},
we will not consider them in the model
formulated here.

Although the full model should
be formulated at large strains as in \cite{RouSte18TEPR},
geophysical applications in solid parts of the Earth (mainly \TTT the crust \EEE
\TTT and the \EEE lithospheric part\TTT{s} \EEE of the mantle) are formulated
at small strains, which can also be more efficiently implemented on computers.
Even, mostly
seismic sources (tectonic earthquakes by fast ruptures of lithospheric faults)
are separated from seismic wave propagation in most of geophysical
simulations, although physically these
two processes are obviously coupled as also captured in the model
presented here. Simultaneously with this small elastic strain assumption
which is well relevant in all processes in the lithosphere,
there might be a large inelastic strain accumulated during slow
tectonic processes on the mentioned large-time scales.
Simultaneously, we will consider inertia so that seismic waves typically
emitted during fast damage and subsequent inelastic shift during
earthquakes are not excluded from the model.

Large inelastic and total strains lead in general also to large displacements.
Then the usual dilemma between \SSS the \EEE Lagrangian and  \SSS the \EEE
Eulerian description arises.
In contrast to \SSS the \EEE standard choice in solid mechanics, we will use
the Eulerian description like suggested essentially in
\cite{LyaBeZ14CDBF,LyHaBZ11NLVE}. Then, all time derivatives in the model
should be convective, i.e.\ the material derivatives. As a consequence, in
particular, the Korteweg-like stresses \FFF arise \EEE from the gradients of internal
variables and the inertial forces
\FFF need \EEE careful formulation and treatment like in fluid mechanics of
so-called quasi-incompressible fluids, cf.\ \cite{Tema69ASEN,Tema77NSE},
refined in the context of elastic ``semi-compressible'' fluids in \TTT the \EEE
consistent Eulerian description in \cite[Sect.5]{Roub20QSF}.

As e.g.\ in \cite{LyaHam07DWFF,LyHaBZ11NLVE,Roub17GMHF}, we use the
Green-Naghdi \cite{GreNag65GTEP}
additive splitting of the total strain but do not
assume the inelastic strain to be small. 
In contrast to \cite{Roub17GMHF}, we formulate the model fully in \TTT the \EEE
Eulerian
setting, so that all time-derivatives are convective (=\,material). Thus,
in contrast to \cite{LyHaBZ11NLVE}, where the structural stress is incomplete
and no energy balance is thus achieved, we have the correct energy balance
rigorously at disposal. An important attribute is that, like in
\cite{LyHaBZ11NLVE}, we formulate the model not in terms of displacements but
rather in terms of velocities and strains.
We admit stored energies which are nonconvex in the elastic
strain like devised in \cite{LyaMya84BECS} to model unstable response of
damaged rocks and used e.g.\ in numerous geophysical articles as e.g.\
\cite{BeZi08CBEF,FMGK13SBFN,LyaBeZ14CDBF,LyHaBZ11NLVE}, and
simultaneously do not 
use a total-strain gradient (which would not be
physically consistent) but only a total-strain-rate gradient.

The goal of this article is to devise the models from
\cite{LyaBeZ14CDBF,LyHaBZ11NLVE} correctly to respect energy
balance and, thus, to allow for rigorous analysis.
Also \cite{Roub17GMHF}, where the inertial term was not formulated
in the convective way and thus the energy balance contained some nonphysical
term and where (rather for analytical reasons but not physically motivated)
the gradient of the total strain was used, \TTT will \EEE thus be improved.
The fluid flow in poroelastic medium, like devised in
\cite{LyaHam07DWFF} without damage gradient,
%structural stresses,
will be consistently incorporated into the %convective
model, too.

\section{The poro-elastodynamical model}\label{sec-2}
%        ~~~~~~~~~~~~~~~~~~~~~~~~~~~~~~~
We consider a continuum whose motion takes place in a fixed region $\Omega$  of space. We denote by $x$ and $t$ the typical point of $\Omega$ and the typical time. The kinematical ingredients (basic variables of the model) are
\begin{align*}
  &\vv&&\text{velocity (valued in $\R^d$)},&&&&&&
  %\\[-.2em]&\Ee,\,\Ep\!\!\!&&\text{elastic and inelastic (plastic-like) strains (valued in $\R^{d\times d}_{\rm sym}$)},&&&&&&
  \\[-.2em]&\Ee&&\text{elastic strain
    (valued in $\R^{d\times d}_{\rm sym}$)},&&&&&&
   \\[-.2em]&\Ep&&\text{inelastic (plastic-like) strain
    (valued in $\R^{d\times d}_{\rm sym}$)},&&&&&&
  \\[-.2em]&\alpha&&\text{other internal variables (as damage, breakage, and/or porosity, valued in $\R^\ell$)},&&&&&&
  \\[-.2em]&\zz&&\text{water (or oil) content (scalar valued)},&&&&&&
 \end{align*}
 with $\R^{d\times d}_{\rm sym}$ denoting the set of symmetric $d{\times}d$-matrices. In addition, we shall use the auxiliary variable $\mu$ which
\SSS will be in a position of a \EEE
%represents
chemical potential\SSS, having here a concrete meaning of the so-called pore
pressure\EEE.

The inelastic strain can incorporate a creep strain to describe Maxwellian rheology or plastic strain to describe activated slip processes which develop, for example, during earthquakes.

We investigate the following system of partial differential 
equations/inclusions:
\begin{subequations}\label{eq}\begin{align}\label{eq:1a}
    &\density\mtd\vv=\operatorname{div}\big(\partial_{\Ee}^{}\varphi(\Ee,\alpha,\zz)
    %{\boldsymbol S}
        +k_{\rm v}\EE(\vv)+{\boldsymbol S}_{\rm str}
    %-\operatorname{div}{\boldsymbol S}_{\rm hpr}
\big)+{\bm f}-
    \frac{\density}2(\operatorname{div}\vv)\vv\,,
%\displaybreak
\\[-.2em]\label{eq:1b}
  &%\bdot
  \mtdEe=\EE(\vv)-\mtd\Ep+k_{\rm e}\Delta\partial_\Ee\varphi(\Ee,\alpha,\zz)  \quad\ \ \text{ with }\ \ \EE(\vv):=\operatorname{sym}\nabla{\bm v}
  %\bdot
\,,
  %\\
      %  &{\bm S}=\varphi'_{\Ee}(\Ee,\alpha)+\Big(\frac \varrho 2 |\bm v|^2+\varphi(\Ee,\alpha,\zz)+\frac {{\color{red}\text{\huge$\kappa$}}_1}2|\nabla\Ee|^2+\frac {{\color{red}\text{\huge$\kappa$}}_2}2|\nabla\Ep|^2+\frac {{\color{red}\text{\huge$\kappa$}}_3}2|\nabla\alpha|^2\Big){\bm I}+\frac {{\color{red}\text{\huge$\kappa$}}_1}2 \nabla\Ee\otimes\nabla..
\\\label{eq:3}
      &\partial_{\mtd{\Ep}}\zeta\Big(\alpha,\zz;\mtdEp,\mtd\alpha\Big)
 -\partial_{\Ee}^{}\varphi(\Ee,\alpha,\zz)
 %{\boldsymbol S}_{_{\rm KV}}
 \ni%k_{\rm v}\Delta\mtdEp
 k_{\rm p}\Delta\Ep\,,
        \\[.1em]\label{eq:4}
        &\partial_{\mtd{\alpha}}\zeta\Big(\alpha,\zz;\mtdEp,\mtd\alpha\Big)+\partial_\alpha^{}\varphi(\Ee,\alpha,\zz)\ni k_{\rm a}\Delta\alpha\,,
 \\[.1em]\label{eq:5}
 &\mtd\zz=\operatorname{div}(\mathbb M(\alpha,\zz)\nabla\mu)
 %\,{\bm j}(\alpha,\zz,\nabla\mu)
 \ \ \ \text{ with }\ \ \ \mu=\partial_{\zz}^{}\varphi(\Ee,\alpha,\zz)\,,
\end{align}\end{subequations}
where we use the conventional notation
\begin{equation}\label{eq:19}
%(\cdot)\!\mtd{^{}}
\mtd{(\cdot)}=\Big[\pdt{}+{\bm v}\cdot\nabla\Big](\cdot)
\end{equation}
to denote the material derivative with respect to time
\INSERT{and where ``$\,\partial\,$'' denotes the partial derivatives or, in
  \FFF(\ref{eq}c,d)\EEE,
  the convex subdifferential to allow for nonsmoothness of the
  dissipation potential $\zeta(\alpha,\zz;\cdot,\cdot)$ at zero rates to model activated processes in inelastic strain and damage/porosity evolution}.
Here $\density>0$ is a reference mass density and $\varphi(\Ee,\alpha,\zz)$ is the free-energy density. Moreover, ${\boldsymbol S}_{\rm str}$ in \eqref{eq:1a}
is the \emph{structural stress} (called also Korteweg's \cite{Kort01FPEM} or Ericksen's \cite{Eric91LCVD} stress) given here as
%\begin{subequations}\label{structural}
  \begin{align}%\nonumber
  {\boldsymbol S}_{\rm str}&=k_{\rm p}\nabla\Ep\boxtimes\nabla\Ep+
     k_{\rm a}\nabla\alpha\boxtimes\nabla\alpha%+{\color{red}\text{\huge$\kappa$}}\nabla\zz\otimes\nabla\zz
%\\&\qquad
     -\Big(\varphi(\Ee,\alpha,\zz)+
  %\frac {k_{\rm e}} 2 |\nabla\Ee|^2+
  \frac{k_{\rm p}}2 |\nabla\Ep|^2
  +\frac{k_{\rm a}}2 |\nabla\alpha|^2
                                        %-\frac\varrho2|{\bm v}|^2
   \Big)\boldsymbol I\,.\label{structural-stress}
  \end{align}%\end{subequations}
  The constants $k_{\rm p}$ and $k_{\rm a}$ appearing in (\ref{eq}c,d)
  and \eqref{structural-stress}
  determine the length-scale of the inelastic strain and of the other internal
  variables.
  In fact, $k_{\rm a}$ can rather be a matrix, expressing different
length-scale for particular internal variables and possible cross-effects.  
\INSERT{The coefficient $k_{\rm v}$ in \eqref{eq:1a} corresponds to the Kelvin-Voigt rheology, but when combined
with a Maxwell rheology which may be governed by \eqref{eq:3},\COMMENT{HERE IT WAS \eqref{eq:4} INCORRECTLY} we actually
obtain the Jeffreys' rheology, as used e.g.\ in \cite{LyHaBZ11NLVE}.}
Moreover,
  $\nabla\Ep\boxtimes\nabla\Ep=\sum_{i,j=1}^d\nabla(\Ep)_{ij}\otimes\nabla(\Ep)_{ij}$,
  i.e.\ component-wise $[\nabla\Ep\boxtimes\nabla\Ep]_{ij}=\sum_{k,l=1}^d\frac{\partial}{\partial x_i}(\Ep)_{kl}\frac{\partial}{\partial x_j}(\Ep)_{kl}$,
  % \COMMENT{OK?? I WOULD ANYHOW USE BETTER $\otimes$ AS A GENERALIZATION OF THE USUAL SCALAR MEANING WITH THIS EXPLANATION. {\color{blue}It is okay, and I would keep the current notation.}}
  and similarly
  $\nabla\alpha\boxtimes\nabla\alpha=\sum_{i=1}^\ell \nabla\alpha_i\otimes\nabla\alpha_i$, whereas ${\boldsymbol I}$ is the identity matrix. 
%  Eventually, in \eqref{eq} we have used a so-called {\it hyper-stress}
%  ${\boldsymbol S}_{\rm hpr}$ as
%% \COMMENT{MAYBE RATHER $k_{\rm e}\nabla\mtd{\Ee}$ SO THAT $\Ee$ NOT OCCURING IN \eqref{structural-stress}?? AND IF $\varphi(\cdot,\alpha)$ convex, then even this not needed ???}
%\begin{align}\label{Shpr}{\boldsymbol S}_{\rm hpr}=k_{\rm v}\nabla\EE(\vv)%\Ee
%  \,\end{align}
%  with $k_{\rm v}$ another viscous capillarity-like modulus
  %  (physically in units Pa\,m$^2$\,s).
  Furthermore, $\zeta(\alpha,\zz;\mtd\Ep,\mtd\alpha)$ in (\ref{eq}c,d)
  is the dissipation potential in general nonsmooth at
  $\mtd\Ep=0$ and $\mtd\alpha=0$.

  The equation \eqref{eq:5} is (the convective variant of) the standard
\INSERT{Fick-type diffusion} driven by the gradient of the chemical potential $\mu$
\TTT with \EEE $\mathbb M=\mathbb M(\alpha,\zz)$ \TTT being \EEE
a positive-definite mobility matrix. The structural force, i.e.\ the last
term in \eqref{eq:1a}, was proposed by
\FFF R.\,Temam \EEE
\cite{Tema69ASEN}, cf.\ also \cite[Ch.\,III,\S\,8]{Tema77NSE}. Beside balancing energetics, this force \FFF vanishes \EEE in the
  incompressible limit, which was the motivation of \cite{Tema69ASEN}. The calculations we perform below provide a justification of this term.
The decomposition
\eqref{eq:1b} is legitimate
\SSS in some special (particular in stratified situations),
cf.\ Remark~\ref{rem-object}, \EEE
below, while the diffusion in \eqref{eq:1b} is discussed in Remark~\ref{rem-dif}.  The structural stresses \eqref{structural-stress} are usually negligible
  but are important, beside balancing energetics,
  ``in narrow zones with high damage gradients or damage fronts separating
  between areas with intact and highly damaged material'', as claimed in
  \cite{LyaBeZ14CDBF}.

We have to complete the system
\eqref{eq}--\eqref{structural-stress}
      %{Shpr}
      by suitable boundary conditions, say
\begin{subequations}\label{BC}
\begin{align}\label{BC1}
  &
  \vv\cdot{\bm n}=0,
  \\\label{BC2}
  &
  %\mathbf P_{\textsc t}
  \big((\partial_{\Ee}^{}\varphi(\Ee,\alpha,\zz){+}k_{\rm v}\EE(\vv){+}{\boldsymbol S}_{\rm str}\TTT)\bm n\EEE
  %-\operatorname{div} \bm S_{\rm hpr}){\bm n}-\divS({\boldsymbol S}_{\rm hpr}{\bm n})
  \big)_{\mathsf t}+\gamma\bm v_{\mathsf t}=\bm g_{\mathsf t},\ \ \ \ \ \
 % \mathbf P_{\textsc t}(
  %{\boldsymbol S}_{\rm hpr}:({\bm n}\otimes{\bm n}) =\bm 0,
  \\\label{BC3}
  &({\bm n}{\cdot}\nabla){\partial_{\Ee}^{}\varphi(\Ee,\alpha,\zz)}=0,\ \ \ \ \
  ({\bm n}{\cdot}\nabla)\Ep=0,\ \ \ \ \ 
 % \\&
  \nabla\alpha{\cdot}{\bm n}=0,\\
  &
  \mathbb M(\alpha,\zz):(\nabla\mu\otimes{\bm n})
  %        {\bm j}(\alpha,\zz,\nabla\mu)\cdot{\bm n}
        %  +{\color{red}\text{\huge$\varkappa$}}\mu={\color{red}\text{\huge$\varkappa$}}\mu_{\rm ext}
=h,
% \qquad\text{ and }\qquad\nabla\zz\cdot{{\bm n}}=0,
    \label{BCe}
\end{align}
\end{subequations}
where
$(\cdot)_{\mathsf t}$ denotes \FFF the \EEE
tangential component of a vector, i.e.\ e.g.\ 
$\bm v_{\mathsf t}=\bm v-(\vv\cdot{\bm n}){\bm n}$
is the tangential velocity (a vector).
Moreover, $\gamma>0$ is a viscous drag coefficient, $\bm g$ a given surface
mechanical load, $h$ is a prescribed inward boundary flux.
  
To unveil the energetic structure of System \eqref{eq}, we derive
the structural \SSS a priori \EEE estimates that we will also use in the analytical
part of this paper. To this effect, we report for the reader's sake some calculations \SSS to be used for the inertial term in \eqref{eq:1a} and the 
higher-order terms in (\ref{eq}c,d)\EEE,
which rely on the representation \eqref{eq:19} of the material time derivative:

\SSS
\begin{lemma}
For $\varrho$ constant and for any sufficiently smooth velocity
field $\vv$ with $\vv{\cdot}\bm n=0$ on the boundary and any smooth
field ${\bm A}$ with $(\bm n{\cdot}\nabla)\bm A=0$ on the boundary, the
following integral identities hold:
\begin{subequations}\begin{align}\label{eq:12}
&\frac{\d}{\d t}\int_\Omega \frac \varrho 2 |\bm v|^2\,\d x
=\int_\Omega\varrho\Big(\mtd{\bm v}+\frac12(\operatorname{div}\bm v){\bm v}\Big)\cdot\bm v\,\d x\quad\text{ and}
\\&\label{eq:14}
\frac{\d}{\d t}\int_\Omega\frac 12 |\nabla{\bm A}|^2\,\d x=
\int_\Omega\Big(\frac 1 2 |\nabla{\bm A}|^2\bm I-\nabla{\bm A}\boxtimes\nabla{\bm A}\Big):\bm e(\bm v)-\Delta{\bm A}:\mtd{\bm A}\,\d x\,.
\end{align}\end{subequations}
\end{lemma}

\EEE

%\begin{proof}
\SSS\noindent{\it Sketch of the proof.} \EEE
The first calculation follows from an application of the divergence theorem
and of the requirement that the normal component of $\bm v$ vanishes on the
boundary, taking into account that the density $\varrho$ is constant:
  \begin{align}%\label{eq:12}
    \int_\Omega\!\!\varrho\mtd{\bm v}{\cdot}\bm v+\frac\varrho 2(\operatorname{div}\bm v)|\bm v|^2\,\d x&= \int_\Omega\!\!\varrho\frac{\partial\bm v}{\partial t}{\cdot}\bm v+\!\!\!\underbrace{\varrho(\bm v{\cdot}\nabla)\bm v{\cdot}\bm v}_{%\displaystyle
      =\bm v{\cdot}\nabla(\varrho|\bm v|^2)/2
      %\frac{\varrho|\bm v|^2}{2}
    }\!\!\!+\operatorname{div}\bm v\frac{\varrho|\bm v|^2\!\!}2\,\d x
\nonumber
%\\[0.6em]&
=\frac{\d}{\d t}\int_\Omega \frac \varrho 2 |\bm v|^2\,\d x
+\int_\Gamma\!\!\!\frac{\varrho|\bm v|^2\!\!}2 \underbrace{\bm v{\cdot}\bm n}_{%\displaystyle
  = 0}\,\d S\,.
  \end{align}
  
  The second %result
  \SSS identity \eqref{eq:14} 
 results from the calculus: \EEE
\begin{align}%\label{eq:14}
    -\int_\Omega &\Delta{\bm A}{:}\mtd{\bm A}\,\d x=\int_\Omega \nabla{\bm A}{:}\nabla\mtd{\bm A}\,\d x-\int_\Gamma \underbrace{(\bm n{\cdot}\nabla){\bm A}}_{%\displaystyle
      =0}{:}\mtd{\bm A}\,\d S
      %\nonumber\\&
  =\int_\Omega \nabla{\bm A}\operatorname{\Vdots}\nabla\frac{\partial {\bm A}}{\partial t}+\nabla{\bm A}\operatorname{\Vdots}\nabla(\bm v{\cdot}\nabla)\bm{\bm A}\,\d x\nonumber\\
    &=\frac{\rm  d}{\d t}\int_\Omega\frac 12 |\nabla{\bm A}|^2\,\d x+\int_\Omega\nabla{\bm A}\boxtimes\nabla{\bm A}:\nabla\bm v+\underbrace{\nabla{\bm A}\operatorname{\Vdots}(\bm v\cdot\nabla)\nabla{\bm A}}_{
      %\displaystyle
      =\bm v\cdot\nabla%\frac{|\nabla{\bm A}|^2}{2}
    |\nabla{\bm A}|^2/2}\,\d x\nonumber\\
    &=\frac{\d}{\d t}\int_\Omega\frac 12 |\nabla{\bm A}|^2\,\d x+\int_\Omega\Big(\nabla{\bm A}\boxtimes\nabla{\bm A}-\frac 1 2 |\nabla{\bm A}|^2\bm I\Big){:}\bm e(\bm v)\,\d x+\int_\Gamma\frac 1 2 |\nabla{\bm A}|^2\underbrace{\bm v\cdot\bm n}_{%\displaystyle
      =0}\,\d S\,.\hspace{5em}\Box\hspace{-2em}
  \nonumber\end{align}
%\end{proof}

\medskip

\TTT Let us remark that \EEE the above result \eqref{eq:12} is consistent
with the interpretation proposed in \cite{Toma21ITST} of Temam's extra force
$\varrho(\operatorname{div}\bm v)\bm v/2$ \cite{Tema69ASEN,Tema77NSE}, an
interpretation based on the requirement that the power expenditure of the
inertial force be equal to minus the rate of change of kinetic energy
\cite{podio1997inertia}.

Taking \TTT now \EEE the scalar product of both sides of \eqref{eq:1a} with the velocity field $\vv$ and integrating over $\Omega$ and using standard divergence identities, as well as \eqref{eq:12}, in combination with the boundary condition \eqref{BC1} and \eqref{BC2} we obtain:
  \begin{align}\label{eq:21}
    0=&\frac{\d}{\d t}\int_\Omega \frac \varrho 2 |\bm v|^2\,\d x+\int_\Omega k_{\rm v}\bm e(\bm v){:}\bm e(\bm v)
    +(\bm S{+}\bm S_{\rm str}){:}\bm e(\bm v)-\bm f{\cdot}\bm v\,\d x
    %\nonumber\\&
    +\int_\Gamma \gamma|\bm v|^2\!-\bm g{\cdot}\bm v\,\d S.
  \end{align}
  On testing \eqref{eq:3} by $\mtd\Ep$ and using, in the order, the boundary conditions for $\bm S=\partial_\Ee(\Ee,\alpha,\chi)$, $\Ep$, and $\bm v$ imposed with \eqref{BC3} and \eqref{BC1}, respectively, making also use of \eqref{eq:14}
  \SSS for ${\bm A}=\Ep$\EEE, we obtain:
  \begin{align}\label{eq:18b}
    0=&\frac{\d}{\d t}\int_{\Omega} \frac{k_{p}}{2}|\nabla \Ep|^{2}\,\d x+\int_\Omega\partial_\Ee\varphi(\Ee,\alpha,\zz):\mtd\Ee+\partial_{\mtd\pi}\zeta\left(\alpha,\chi;\mtd\Ep,\mtd\alpha\right):\mtd\Ep+k_{\rm e}\nabla\bm S\Vdots\nabla\bm S\nonumber\\
    &\hspace{14em}+\Big(k_{\rm p}\nabla\Ep\boxtimes\nabla\Ep-\frac {k_{\rm p}}2|\nabla\Ep|^2\bm I\Big):\bm e(\bm v) -\bm S:\bm e(\bm v)\,\d x.
  \end{align}  
  In a similar fashion, we operate on \eqref{eq:4},\COMMENT{HERE IT WAS \eqref{eq:3} INCORRECTLY} using $\mtd\alpha$ as test function \SSS and \eqref{eq:14}
  for ${\bm A}=\alpha$\EEE, and taking into account the Neumann boundary condition for $\alpha$ in \eqref{BC3}. The result is:
           \begin{align}\label{eq:17}
    0=&\frac{\d}{\d t}\int_{\Omega} \frac{k_{a}}{2}|\nabla \alpha|^{2}\,\d x+\int_\Omega\partial_\alpha\varphi(\Ee,\alpha,\zz):\mtd\alpha+\partial_{\mtd\alpha}\zeta\left(\alpha,\chi;\mtd\Ep,\mtd\alpha\right)\cdot\mtd\alpha\nonumber\\
    &\hspace{14em}+\Big(k_{\rm a}\nabla\alpha\boxtimes\nabla\alpha-\frac {k_{\rm a}}2|\nabla\alpha|^2\bm I\Big):\bm e(\bm v)\,\d x.
           \end{align}
           Finally, on testing the equations in \eqref{eq:5}, respectively, by $\mu$ and $\mtd\zz$, and adding the resulting equations, we obtain
\begin{align}
   0=\int_\Omega \partial_\zz\varphi(\Ee,\alpha,\zz)\mtd\zz
   +\mathbb M\nabla\mu:\nabla\mu
   % +\Big({\color{red}\text{\huge$\kappa$}}\nabla\zz\boxtimes\nabla\zz-\frac {\color{red}\text{\huge$\kappa$}} 2 |\nabla\zz|^2\bm I\Big):\bm e(\bm v)
  \,\d x
                 %\underbrace{-\int_\Gamma\mu\mathbb M:\nabla\mu\otimes\bm n\,\d S}_{\displaystyle=
                 -\int_\Gamma h\mu\,\d S\,,
%}
  \label{eq:15b}
\end{align}
where we used also the boundary condition \eqref{BCe}. By summing the
estimates %\eqref{eq:18b}, \eqref{eq:17}, and \eqref{eq:15b},
\eqref{eq:18b}--\eqref{eq:15b} and by observing that
  \begin{align}\label{eq:20}
    &\int_\Omega \partial_\Ee\varphi(\Ee,\alpha,\zz):\mtd\Ee
    +\partial_\alpha\TTT\varphi\EEE(\Ee,\alpha,\zz):\mtd\alpha+\partial_\zz\TTT\varphi\EEE(\Ee,\alpha,\zz):\mtd\zz\,\d x\nonumber\\
    &=\int_\Omega \frac{\partial\varphi_\Ee(\Ee,\alpha,\zz)}{\partial t}+
    %\int_\Omega
    (\bm v\cdot\nabla(\varphi_\Ee(\Ee,\alpha,\zz))\,\d x\nonumber\\
    &=\frac\d{\d t}\int_\Omega \frac{\partial\varphi_\Ee(\Ee,\alpha,\zz)}{\partial t}\,\d x-\int_\Omega \varphi_\Ee(\Ee,\alpha,\zz)\bm I:\bm e(\bm v)\,\d x+\int_\Gamma\varphi_\Ee(\Ee,\alpha,\zz)\underbrace{\bm v\cdot\bm n}_{\displaystyle =0}\,\d S\,,
  \end{align}
  we obtain the partial energy balance
    \begin{align}\label{eq:18}
      &\frac{\d}{\d t}\int_\Omega\varphi(\Ee,\alpha,\zz)+\frac{k_{p}}{2}|\nabla \Ep|^{2}+\frac{k_{\rm a}}{2}|\nabla\alpha|^2\,\d x
        \nonumber \\&+\int_\Omega\partial_{\mtd\pi}\zeta\left(\alpha,\chi;\mtd\Ep,\mtd\alpha\right):\mtd\Ep+\partial_{\alpha}\zeta\left(\alpha,\chi;\mtd\Ep,\mtd\alpha\right):\mtd\alpha+k_{\rm e}\nabla\bm S\Vdots\nabla\bm S-(\bm S+\bm S_{\rm str}):\bm e(\bm v)\,\d x\nonumber\\
      &\quad=\int_\Gamma
       h\mu\,\d S.
    \end{align}
On summing \eqref{eq:18} and \eqref{eq:21}, the contributions from the thermodynamic stress $\bm S$ and the structural stress $\bm S_{\rm str}$ cancel, and we arrive at the following total energy balance:
      \begin{align}\nonumber
&\!\!\int_\Omega\!\linesunder{\frac\varrho2|\vv(t)|^2}{kinetic}{energy}\!\!+
  \!\!\lineunder{\varphi(\Ee(t),\alpha(t),\zz(t))
   +\frac{k_{\rm p}}2|\nabla\Ep|^2+\frac{k_{\rm a}}2|\nabla\alpha(t)|^2
  }{stored energy}\!\d x
%\\&\nonumber\qquad
  \\&\qquad\nonumber
  +\int_0^t\!\!\bigg(\int_\Omega\!\!\lineunder{\xi\Big(%\Ee,
    \alpha,\zz;
  %\mtdEe,
    \mtdEp,%\EE(\vv),
    \mtd\alpha%,\nabla\mu,\nabla{\bm S}
    %%\Ee,\nabla\alpha,\nabla\zz
    \Big)
    %+\bbM(\zz)\nabla\mu\cdot\nabla\mu
+k_{\rm v}|\EE(\vv)|^2
+{\mathbb M}(\alpha,\zz)\nabla\mu{\cdot}\nabla\mu
+k_{\rm e}\big|\nabla{\bm S}\big|^2
  }{bulk dissipation rate}% and}{diffusion of the phase field}
  \!\!\d x
  +\int_\Gamma\!\!\!\!\!\!\!\!\!\linesunder{\gamma|\vv_{\rm t}|^2}{boundary}{dissipation rate}\!\!\!\!\!\!\!\!\!\d S\bigg)\d t
\\[-.1em]&\ \nonumber
 =\int_0^t\!\!\bigg(\int_\Omega\!\!\!\!\!\linesunder{{\bm f}{\cdot}\vv_{_{_{}}}}{power of}{bulk load}\!\!\!\!\!\d x
 +\int_\Gamma\!\!\!\linesunder{{\bm g}_{\rm t}{\cdot}\vv_{\rm t}
   +h\mu_{_{_{}}}}{power of}{boundary load}\!\!\!\!\d S\bigg)\,\d t
\\[-.7em]&\qquad\qquad\qquad\qquad\qquad\qquad\qquad\qquad\ \
  +\int_\Omega\frac\varrho2|\vv_0|^2+\varphi({\Ee}_0,\alpha_0,\zz_0)
   \SSS+\frac{k_{\rm p}}2|\nabla\Ep|^2\EEE
   +\frac {k_{\rm a}} 2 |\nabla\alpha_0|^2\,\d x\,,
  \label{energy}\end{align}
  where we abbreviated
 \begin{align}\label{def-of-xi}
   \xi(\alpha,\zz;\mtd{\Ep},\mtd\alpha)
   =\partial_{\mtd\Ep}\zeta\Big(\alpha,\zz;\mtd{\Ep},\mtd\alpha\Big):\mtd{\Ep}
   %+k_{\rm v}\nabla\mtdEp:\nabla\mtdEp
   +\partial_{\mtd\alpha} \zeta\Big(\alpha,\zz;
   \mtd{\Ep},\mtd\alpha\Big)\cdot\mtd\alpha\,.
   \end{align}

%\section{Remarks}
 %        ~~~~~~~

 \def\OMEGA{\text{$\vspace*{-.2em}\includegraphics[width=1.em]{omega.pdf}$}}

\SSS
\begin{remark}[Additive elasto-inelastic strain rate decomposition]
The concept of additive strain-rate of the type like \eqref{eq:1b}
%decomposition of the type $\EE(\vv)=\mtd{}\Ee+\mtd{}\Ep$, cf.\ \eqref{eq:1b},
dates back basically
to Hill \cite{Hill58GTUS} and Prager \cite{Prag61EDDS}.  This 
concept has been widely used in literature, cf.\ e.g.\
\cite[Sec.8.6]{HasYam13IFST} or \cite[Sec.8.3]{Maug92TPF}. For discussing
some limitations, see e.g.\ \cite{JiaFis17ADRD}. It should however be
noticed that our equation \eqref{eq:1b} only mimicks the standard additive
decomposition in a simplified way that we write such decomposition like in the small strain setting. 
\end{remark}
\EEE

\begin{remark}[Objective variant.]\label{rem-object}
\SSS There is another simplification we used in the decomposition
\eqref{eq:1b}. \EEE It is known that the material derivatives of tensors $\Ee$ and
$\Ep$ used in (\ref{eq}b-d) are not frame indifferent\SSS, and that \EEE
 % Unless $\bm F\sim\bm I$
 % as in Remark~\ref{rem-decomp},
one should rather use Oldroyd's or Zaremba-Jaumann's \TTT or Green-Naghdi's \EEE
time derivative, cf.\ e.g.\
%\cite[Sect.5.5]{Mart19PCM}
\cite{Haup02CMTM,HasYam13IFST,Mart19PCM}. See also
\cite[Sect.\,5.4]{surana}. \SSS Nevertheless, the simplifying purely convective
but non-objective variant used here is also often exploited in geophysical
modelling, although the corresponding symmetric structural stress is usually
not reflected correctly there, cf.\ e.g.\
\cite{BilHir07RCSD,LyHaBZ11NLVE,ReLYue03MSZG}. A certain legitimacy of this
simplification is in stratified simple-shear situations, as articulated in
\cite[Proposition~1]{Roub17GMHF}. \SSS There it is shown, starting from the
Kr\"oner-Lee multiplicative decomposition \cite{Kron60AKVE,LeeLiu67FSEP},
that if displacements are large in one direction (parallel to the
stratification) and if the plastic part is a simple shear, then the additive
decomposition holds up to higher-order terms. Cf.\ also
\cite[Sect.\,8.1.3]{Haup02CMTM} or \cite[Sect.\,5]{Volo13AELD}. {\EEE} Actually,
the objective variant (most suitably using Zaremba-Jaumann
derivatives) is amenable for analysis, too;
%when the corresponding structural stress would be nonsymmetric and the above calculus would be enhanced correspondingly,
cf.\ \cite{Roub??SPTC,Roub??TCMP}.  
\EEE

  %\COMMENT{HERE MAYBE DETAILS - Jaumann-Zaremba BETTER (because for
 %    fluids it gives just the convective derivative of the pressure) MAYBE ALSO RELEVANT:\\
%     Z. Fiala: Geometrical setting of solid mechanics. Annals of Physics 326 (2011) 1983--1997
%     \\OR\\Z. Fiala: Geometry of finite deformations and time-incremental analysis. International Journal of Non-Linear Mechanics 81 (2016) 230--244\\OR\\Z. Fiala: Objective time derivatives revised. Z. Angew. Math. Phys. (2020) 71:4
%     \\
%     FOR Oldroyd e.g.:
%   Xianpeng Hu and Hao Wu:   Long-time behavior and weak-strong uniqueness for
% incompressible viscoelastic flows.  arXiv:1411.0518, 2014
%   }
\end{remark}
  
\begin{remark}[The stress diffusion in \eqref{eq:1b}.]\label{rem-dif}
  The $k_{\rm e}\Delta\partial_\Ee\varphi$-term in the ``geometrical'' equation
  might be rather controversial. In fluid dynamics, this stress diffusion was
  advocated in series of
works by H.\,Brenner, cf. e.g.\ \cite{Bren05KVT,Bren06FMR}. Independently, such
  diffusion was used also in \cite{BFLS18EWSE,LiLiZh05HVF}. % p. 1461 in LiLiZh05HVF  
  Cf.\ also a discussion in \cite{OtStLi09IDCM} and a  thermodynamical
  justification in \cite{VaPaGr17EMFF}. It is also somehow similar to the recent proposal in \cite{burczak2017}. When omitting information about
  $\nabla\partial_\Ee\varphi$ and thus about $\nabla\Ee$, our \SSS a priori \EEE estimates would work even with $k_{\rm e}=0$. A
  positive (even arbitrarily small) $k_{\rm e}>0$ together with the convexity of
$\varphi(\cdot,\alpha,\zz)$ is needed only for facilitating
convergence of approximate solutions. Therefore, this questionable regularizing
diffusion will not essentially influence global energetics and presumably
would not be seen on numerically stable algorithms if $k_{\rm e}>0$ would be
small.
\end{remark}

\section{Analysis of an initial-value problem for \eqref{eq}--\eqref{BC}}\label{sec-3}
%       ~~~~~~~~~~~~~~~~~~~~~~~~~~~~~~~~~~~~~~~~~~~~~~~~~~~~~~~~~~~~~~~

%\COMMENT{%SHALL WE PUT $d=3$?
%A GENERAL $d\le3$ CONSIDERED}
We carry out the analysis for $d\le 3$. We will now consider the evolution governed by the system \eqref{eq}
on a fixed time interval $I=[0,T]$ and complete \eqref{eq} by initial conditions
\begin{align}\label{IC}
\vv|_{t=0}^{}=\vv_0,\ \ \ \Ee|_{t=0}^{}={\Ee}_0,\ \ \ 
  \Ep|_{t=0}^{}={\Ep}_0,\ \ \ 
\alpha|_{t=0}^{}=\alpha_0,\ \ \ \zz|_{t=0}^{}=\zz_0\,.
\end{align}

We will perform the analysis by a rather constructive approximation,
which can in principle be also used as a conceptual numerical algorithm
for which we will prove numerical stability and convergence at least in
terms of subsequences of approximate solutions.

We will use the standard notation concerning the Lebesgue and the Sobolev
spaces, namely $L^p(\varOmega;\R^n)$ for Lebesgue measurable functions $\varOmega\to\R^n$ whose Euclidean norm is integrable with $p$-power, and
$W^{k,p}(\varOmega;\R^n)$ for functions from $L^p(\varOmega;\R^n)$ whose
all derivative up to the order $k$ have their Euclidean norm integrable with
$p$-power. Moreover, for a Banach space $X$
and for $I=[0,T]$,
we will use the notation $L^p(I;X)$ for the Bochner space of Bochner
measurable functions $I\to X$ whose norm is in $L^p(I)$, 
and $W^{1,p}(I;X)$ for functions $I\to X$ whose distributional derivative
is in $L^p(I;X)$. Furthermore, $C_{\rm w}(I;X)$ will denote the
Banach space of weakly continuous functions $I\to X$.
We also write briefly $H^k=W^{k,2}$. 

We will assume, with some $\epsilon>0$ arbitrarily small, that
%\COMMENT{\eqref{ass:1-} MAYBE FOLLOWS NOW FROM \eqref{ass:1} EVEN WITH LESSER EXPONENT ???}
\begin{subequations}\label{ass}\begin{align}\nonumber
    &\varphi:\R^{d\times d}\times\R^\ell\times\R\to\R\ \text{ twice continuously
      differentiable,
    %}\\&\nonumber\qquad\qquad\qquad\text{separately convex,
      bounded from below with}
    \\&\label{ass:1-}\qquad
    |\partial\varphi(\Ee,\alpha,\zz)|\le (1+|\Ee|^{3/2-\epsilon}\!
    +|\alpha|^{3-\epsilon}\!+|\zz|^{3-\epsilon})/\epsilon\ \ \text{ and }\ \
    \\&\label{semi-convex}\qquad
    %\exists\ell>0:\ \ \ 
 (\Ee,\alpha,\zz)\mapsto\varphi(\Ee,\alpha,\zz)+\frac1{2\epsilon}|\alpha|^2
    %\frac12K(|\Ee|^2+|\alpha|^2+|\zz|^2)
    \ \text{ is convex},
    \\&\nonumber\qquad
      \partial_{\Ee\alpha}^2\varphi,\ \partial_{\Ee\zz}^2\varphi\ \text{ bounded and }\ \varphi(\cdot,\alpha,\cdot)\ \text{ uniformly convex, i.e.}\ 
    %\\&\nonumber\qquad
    \forall\Ee,\widetilde{\bm E}\in\R_{\rm sym}^{d\times d}, \ \alpha\in\R^\ell,\ \zz,\widetilde\zz\in\R:\
\\&\qquad\quad\label{ass:1}
    \partial_{\Ee\Ee}^2\varphi(\Ee,\alpha,\zz)\widetilde{\bm E}:\widetilde{\bm E}+2\partial_{\Ee\zz}^2\varphi(\Ee,\alpha,\zz)\widetilde{\bm E}\widetilde{\zz}+\partial_{\zz\zz}^2\varphi(\Ee,\alpha,\zz)\widetilde{\zz}^2%\widetilde{\zz}
%\nonumber\\&\qquad\quad
    \ge\epsilon\big(|\widetilde{\bm E}|^2+|\widetilde{\zz}|^2\big)\,,
\\&\nonumber
\zeta:\R^\ell\times\R\times\R^{d\times d}\times\R^\ell\to\R\
\text{ continuous with}
%\\&\nonumber\qquad
%\zeta\big(\alpha,\zz;\DT{\Ep},\DT\alpha\big)=
%\zeta_{\alpha,\zz}^{(1)}\big(\DT{\Ep}\big)+\zeta_{\alpha,\zz}^{(2)}\big(\DT\alpha\big)
%%\zeta(\alpha,\zz;\cdot,\cdot):(\R^{d\times d}{\setminus}\{0\})\times(\R^\ell{\setminus}\{0\})\to\R
%\ \text{ with $\ \zeta_{\alpha,\zz}^{(1)}:\R^{d\times d}{\setminus}\{0\}\to\R$}
%\\&\nonumber\hspace{17.2em}
%  \text{ and $\ \ \zeta_{\alpha,\zz}^{(2)}:\R^\ell{\setminus}\{0\}\to\R$
%  continuously differentiable},
\\&\nonumber\qquad
 \zeta(\alpha,\zz;\cdot,\cdot):\R^{d\times d}\times\R^\ell\to\R
\ \text{ convex and}
\\&\qquad\TTT
\zeta(\alpha,\zz;\cdot,\DT\alpha):\R^{d\times d}\setminus\{0\}\to\R\ \ \ \text{ is
 continuously differentiable}\,,\EEE
\label{smoothnes-1}
\\&\qquad\TTT
\zeta(\alpha,\zz;\DT\Ep,\cdot):\R^\ell\setminus\{0\}\to\R\ \ \ \ \ \text{ is
 continuously differentiable}\,,\EEE
\label{smoothnes-2}
\\&\qquad\epsilon\big(\big|\DT\Ep\big|^2+\big|\DT\alpha\big|^2\big)\le
\zeta\big(\alpha,\zz;\DT\Ep,\DT\alpha\big)\le
\big(1+\big|\DT\Ep\big|^2+\big|\DT\alpha\big|^2\big)/\epsilon\,,
%\ \ \ \text{ with }\ 1<r\le 4/3,
\label{ass:2}
\\\label{ass:3}&
%{\bm j}:\R^\ell\times\R\times\R^d\to\R^d\ \text{ continuous, .....},\ \
%{\bm j}(\alpha,\zz,{\bm m})\cdot{\bm m}\ge\epsilon|{\bm m}|^4\,,
\mathbb M:\R^\ell\times\R\to\R^{d\times d}
\ \text{ continuous, bounded, uniformly positive definite},
\\\label{ass:4}
&\varrho,\ %\text{\COMMENT{MAYBE $\varrho=\varrho(x)$ WOULD BE NICER }}
k_{\rm v},\ k_{\rm p},\  k_{\rm a},\ k_{\rm e}>0\,,
\\\label{ass:5}
&\vv_0\!\in\! L^2(\Omega;\R^d),\ \
        %{\Ee}_0\!\in\! L^2(\Omega;\R^{d\times d}),\
  {\Ep}_0\!\in\! H^1(\Omega;\R_{\rm sym}^{d\times d}),\ \
\alpha_0\!\in\! H^1(\Omega;\R^\ell),\ \ \zz_0\!\in\! L^2(\Omega),
\\\label{ass:6}&{\bm f}\in L^1(I;L^2(\Omega;\R^d))\,,\ \ \ \ \
{\bm g}\in L^2(I;L^\infty(\Gamma;\R^d))\,,\ \ \ \ h\in L^2(I;L^{6/5}(\Gamma))\COMMENT{To check}\,.
\end{align}\end{subequations}
We have formulated our growth assumption \eqref{ass:1} to be valid for
$d=3$ and $d=2$ too, but for the latter case it can be weakened.
%\COMMENT{TO OMIT IF $d=3$ WOULD BE FIXED}
Let us emphasize that we do not assume $\varphi$ convex, which
allows to treat real damage model where $\varphi$ is always nonconvex,
cf. also Remark~\ref{rem-example-energy}.

%Let us note that \eqref{ass:3} does not admit a flux linearly dependent
%on ${\bm m}=\nabla\mu$; this is needed for making the structural stress
%$L^2$ in time and thus making $\pdt\vv$ in duality with $\vv$.

\begin{definition}[Weak solutions to \eqref{eq}--\eqref{BC} with
    \eqref{IC}.]\label{def}
  The 5-tuple $(\vv,\Ee,\Ep,\alpha,\zz)$ with
\begin{subequations}\label{cont-of-sln}\begin{align}
    &\vv\in C_{\rm w}(I;L^2(\Omega;\R^d))\cap L^2(I;H^1(\Omega;\R^d))\ \ \
    \text{\TTT with $\bm n{\cdot}\vv=0$ on $I{\times}\Gamma$\EEE}\,,
    \\&\Ee\in C_{\rm w}(I;L^2(\Omega;\R_{\rm sym}^{d\times d}))
    \,\cap\,L^2(I,H^1(\Omega;\R_{\rm sym}^{d\times d}))\,,
  \\\label{Ep-weak-cont}
  &\Ep\in C_{\rm w}(I;H^1(\Omega;\R_{\rm sym}^{d\times d}))\cap W^{1,4/3}(I;L^2(\Omega;\R_{\rm sym}^{d\times d}))\,,
 \\&\alpha\in C_{\rm w}(I;H^1(\Omega;\R^\ell))\cap W^{1,4/3}(I;L^2(\Omega;\R^\ell))\,,
 \\&\zz\in C_{\rm w}(I;L^2(\Omega))\,,\ \text{ and }\ \mu=\partial_\zz(\Ee,\alpha,\zz)\in L^2(I;H^1(\Omega))
\end{align}\end{subequations}
 will be called a weak solution to the boundary-value problem
 \eqref{eq}--\eqref{BC} with the initial conditions \eqref{IC} if
 ${\bm S}=\partial_\Ee^{}\varphi(\Ee,\alpha,\zz)\in L^2(I;H^1(\Omega;\R^{d\times d}))$,
% $\Delta\Ep\in L^2(I{\times}\Omega;\R^{d\times d})$,
 $\partial_\alpha^{}\varphi(\Ee,\alpha,\zz)\in L^2(I{\times}\Omega;\R^\ell)$,
 %$\Delta\alpha\in L^2(I{\times}\Omega;\R^\ell)$,
 and \FFF the following four integral identities hold: \EEE %if 
  \begin{subequations}\label{def+}\begin{align}\nonumber
      &\int_0^T\!\!\!\int_\Omega\density\Big(\!(\vv{\cdot}\nabla)\vv
      {+}\frac12({\rm div}\,\vv)\vv\Big){\cdot}\widetilde\vv
      +\big({\bm S}
      %\partial_{\Ee}^{}\varphi(\Ee,\alpha,\zz)
      {+}{\boldsymbol S}_{\rm str}\big){:}\EE(\widetilde\vv)+%{\boldsymbol S}_{\rm hpr}
k_{\rm v}\EE(\vv):\EE(\widetilde\vv)
      -\density\vv{\cdot}\pdt{\widetilde\vv}\,\d x\d t
      \\&\hspace{14em}=\int_\Omega\density\vv_0\cdot\widetilde\vv(0)\,\d x
      +\int_0^T\!\!\!\int_\Omega{\bm f}\cdot\widetilde\vv\,\d x\d t
      +\int_0^T\!\!\!\int_\Gamma{\bm g}_{\rm t}\cdot\widetilde\vv_{\rm t}\,\d S\d t
      \intertext{with ${\boldsymbol S}_{\rm str}$ %and ${\boldsymbol S}_{\rm hpr}$
        from \eqref{structural-stress} %--\eqref{Shpr}
        for all $\widetilde\vv\in H^1(I{\times}\Omega;\R^d)
        %;L^2(\Omega;\R^d))\cap L^2(I;H^2(\Omega;\R^d))
        $ with \TTT $\bm n{\cdot}\widetilde\vv=0$ on $I{\times}\Gamma$ and with \EEE $\widetilde\vv(T)=0$, %and if
        }
&\int_0^T\!\!\!\int_\Omega
      (\vv{\cdot}\nabla)(\Ee{+}\Ep):\widetilde{\bm E}
      +k_{\rm e}\nabla{\bm S}\Vdots\nabla\widetilde{\bm E}
      -(\Ee{+}\Ep):\pdt{\widetilde{\bm E}}
%      (\vv\otimes\Ee)\Vdots\nabla\widetilde{\bm E}+\Big(({\rm div}\,\vv)\Ee
%-(\vv{\cdot}\nabla)\Ep-\EE(\vv)\Big):{\widetilde{\bm E}}
\,\d x\d t=\int_\Omega(\Ee_0{+}\Ep_0):\widetilde{\bm E}(0)\,\d x
      \label{eq:1b-weak}
      \intertext{holds for all $\widetilde{\bm E}\in H^1(I{\times}\Omega;\R^{d\times d})$ with $\widetilde{\bm E}(T,\cdot)=0$, %and if
%\COMMENT{PROBLEM:\ CONVERGENCE IN $\partial_\alpha^{}\varphi(\Ee,\alpha,\zz)\cdot A$ IF $\Ee$ ONLY WEAKLY!}
      }
     &\nonumber\int_0^T\!\!\!\int_\Omega
      \zeta(\alpha,\zz;\widetildeEp,\widetilde{\alpha})-{\bm S}
      %\partial_{\Ee}^{}\varphi(\Ee,\alpha,\zz)
      :\Big(\widetildeEp{-}\mtdEp\Big)+\partial_\alpha^{}\varphi(\Ee,\alpha,\zz)
      \cdot\Big(\widetilde{\alpha}{-}\mtd\alpha\Big)
      + k_{\rm p}\nabla\Ep\Vdots\nabla\widetildeEp
      \\[-.3em]&\nonumber\qquad
     % +k_{\rm p}\nabla\Ep\Vdots\nabla\Big(\widetildeEp{-}\mtdEp\Big)
     % +k_{\rm a}\nabla\alpha:\nabla\Big(\widetilde{\alpha}{-}\mtd\alpha\Big)
%-k_{\rm p}\Delta\Ep:\big(\widetildeEp{-}(\vv{\cdot}\nabla)\Ep\big)
      %-k_{\rm a}\Delta\alpha\cdot\big(\widetilde{\alpha}{-}(\vv{\cdot}\nabla)\alpha\big)
+ k_{\rm p}\Delta\Ep{:}(\vv{\cdot}\nabla)\Ep
      +k_{\rm a}\nabla\alpha{:}\nabla\widetilde{\alpha}
      +k_{\rm a}\Delta\alpha{\cdot}(\vv{\cdot}\nabla)\alpha\big)
      \,\d x\d t+\int_\Omega\frac{k_{\rm p}}2|\nabla\Ep_0|^2+\frac{\SSS k_{\rm a}\EEE}2|\nabla\alpha_0|^2\,\d x
    \\[-.1em]&\qquad\qquad\qquad\qquad\qquad\ge
\int_\Omega\frac{k_{\rm p}}2|\nabla\Ep(T)|^2+\frac{\SSS k_{\rm a}\EEE}2|\nabla\alpha(T)|^2\,\d x
+\int_0^T\!\!\!\int_\Omega\zeta\Big(\alpha,\zz;\mtdEp,\mtd\alpha\Big)\,\d x\d t
   \label{VI}
%       &\nonumber\int_0^T\!\!\!\int_\Omega
%          \zeta(\alpha,\zz;\varPi,\widetilde{\alpha})-\partial_{\Ee}^{}\varphi(\Ee,\alpha,\zz):\varPi+\partial_\alpha^{}\varphi(\Ee,\alpha,\zz)\cdot \widetilde{\alpha}
%         +k_{\rm p}\nabla\Ep\Vdots\nabla\varPi
%          +k_{\rm a}\nabla\alpha:\nabla \widetilde{\alpha}\,\d x\d t
%  \\[-.3em]&\nonumber\hspace*{5em}+\int_\Omega\varphi({\Ee}_0,\alpha_0,\zz_0)
%  +\frac{k_{\rm p}}2|\nabla{\Ep}_0|^2+\frac{k_{\rm a}}2|\nabla\alpha_0|^2
%  +\frac{\color{red}\text{\huge$\kappa$}}2|\nabla\zz_0|^2 +\frac{\density}2|\vv|^2\,\d x
%\\[-.3em]&\nonumber\ge\int_\Omega\varphi(\Ee(T),\alpha(T),\zz(T))
%  +\frac{k_{\rm p}}2|\nabla\Ep(T)|^2+\frac{k_{\rm a}}2|\nabla\alpha(T)|^2
% +\frac{\color{red}\text{\huge$\kappa$}}2|\nabla\zz(T)|^2 +\frac{\density}2|\vv(T)|^2\,\d x
%\\[-.3em]&\hspace*{5em}+\int_0^T\!\!\!\int_\Omega\zeta(\alpha,\zz;\mtdEp,\mtd\alpha)+k_{\rm v}|\nabla\EE(\vv)|^2+\mathbb M(\alpha,\zz)\nabla\mu\cdot\nabla\mu\,\d x\d t+\int_0^T\!\!\!\int_\Gamma...................\,\d S\d t\label{VI+}
\intertext{holds for all
$(\widetildeEp,\widetilde{\alpha})\in L^2(I;H^1(\Omega;\R^{d\times d}{\times}\R^{\ell}))$, and %if moreover
}
&
\int_0^T\!\!\!\int_\Omega\big(\mathbb M(\alpha,\zz)\nabla\mu-\zz\bm v\big)\cdot\nabla z
-\zz\pdt z-(\operatorname{div}\vv)\zz z\,\d x
=\int_0^T\!\!\!\int_\Gamma h z\,\d S\d t+\int_\Omega\zz_0z\,\d x
\end{align}\end{subequations}
holds for all $z\in C^1(I{\times}\Omega)$ with $z|_{t=T}^{}=0$ and with
$\mu=\partial_{\zz}^{}\varphi(\Ee,\alpha,\zz)$ a.e.\ in $I{\times}\Omega$.
\end{definition}

\medskip

Let us note that, for the inequality \eqref{VI}, we used the standard
definition of the convex subdifferential of $\zeta(\alpha,\zz;\cdot,\cdot)$
combined with the calculus
\begin{align*}
\int_0^T\!\!\!\int_\Omega\Delta\Ep{:}\Big(\widetildeEp{-}\mtd\Ep\Big)\,\d x\d t
=\frac12\int_\Omega|\SSS\nabla\EEE\Ep(T)|^2-|\SSS\nabla\EEE\Ep(0)|^2\,\d x
-\int_0^T\!\!\!\int_\Omega\nabla\Ep\Vdots\nabla\widetildeEp
+\Delta\Ep{:}(\vv{\cdot}\nabla)\Ep\,\d x\d t\,;
\end{align*}
\SSS{}for the analytical legitimacy of this formula if
$\Delta\Ep$ and $\mtd\Ep$ belong to $L^2(I{\times}\Omega;\R^{d\times d})$
see e.g.\ \cite[Formula (12.133b)]{Roub13NPDE}. \EEE
An analogous calculus for $\alpha$, which both will be actually
legitimate when showing that both $\Delta\TTT\Ep\EEE$ and $\mtd{\TTT\Ep\EEE}$
belong to $L^2(I{\times}\Omega;\R^{d\times d})$ and similarly both $\Delta\alpha$
and $\mtd\alpha$ belong to $L^2(I{\times}\Omega;\R^\ell)$.
\SSS The inequality in \eqref{VI} is also well consistent with the
weak continuity in (\ref{cont-of-sln}c,d) and thus weak lower semicontinuity
of the right-hand side of \eqref{VI}. \EEE

%For the energy balance \eqref{energy} below, let us still define the local
%dissipation rate $\xi$ which is revealed from \eqref{def-of-D} written, by
%using the Green formula, as
%\begin{align}\nonumber \xi(\Ee,\alpha,\zz;
  %%\mtdEe,
  %\mtd\Ep,\EE(\vv),\mtd\alpha,\nabla\mu,\nabla{\bm S}%\Ee,\nabla\alpha,\nabla\zz
  %)&=\partial_{\mtd\Ep}\zeta\big(\alpha,\zz;
  %%\mtdEe,
  %\mtd\Ep,\mtd\alpha\big){:}\mtd\Ep+\partial_{\mtd\alpha} \zeta\big(\alpha,\zz;
  %%\mtdEe,
  %\mtd\Ep,\mtd\alpha\big){\cdot}\mtd\alpha +k_{\rm v}|\nabla\EE(\vv)|^2
%\\&%+k_{\rm v}|\nabla\mtd\Ep|^2
%\ +{\mathbb M}(\alpha,\zz)\nabla\mu{\cdot}\nabla\mu
%+k_{\rm e}\big|\nabla{\bm S}\big|^2
%%\Vdots\nabla\Ee\,.
%%J(\alpha,\zz,\nabla\mu)
%\ \ \text{ with }\ {\bm S}=\partial_{\Ee}^{}\varphi(\Ee,\alpha,\zz)\,.
%\label{def-of-D+}\end{align}
%%with $J$ from \eqref{def-of-J}.
%Let us note that the last term involves
%\begin{align}\nonumber
%\nabla{\bm S}=\nabla\partial_{\Ee}^{}\varphi(\Ee,\alpha,\zz)&=
%\partial_{\Ee\Ee}^2\varphi(\Ee,\alpha,\zz)\nabla\Ee
%\\&\ \ \ +\partial_{\Ee\alpha}^2\varphi(\Ee,\alpha,\zz)\nabla\alpha
%+\partial_{\Ee\zz}^2\varphi(\Ee,\alpha,\zz)\nabla\zz\,,
%\label{calculus-for-diffusion-Ee}\end{align}
%which reveals the dependence on $\nabla\alpha$, and $\nabla\zz$ which is
%``optically'' hidden in the right-hand side of \eqref{def-of-D+}.

\medskip

\begin{theorem}[Existence of weak solutions]
Let $\Omega\subset\R^d$ be Lipschitz and the assumptions \eqref{ass} hold.
Then there exists at least one weak solution $(\vv,\Ee,\Ep,\alpha,\zz)$
to the initial-boundary-value problem \eqref{eq}--\eqref{BC} with \eqref{IC}
according the Definition~\ref{def} which, moreover, satisfies also
$\Delta\Ep\in L^2(I{\times}\Omega;\R^{d\times d})$,
$\Delta\alpha\in L^2(I{\times}\Omega;\R^\ell)$, and
$\nabla\zz\in L^2(I{\times}\Omega;\R^d)$.
\end{theorem}

%\begin{proof}[Sketch of the proof]
\noindent{\it Sketch of the proof.}
For clarity, we divide the proof into {four} steps.
 % First, in Steps 1--4, we will assume $\varphi$ convex and only in Step 5 we
 % use the full generality for a semiconvex $\varphi$.

  \medskip\noindent{\it Step 1. (Approximate solutions - existence)}:
  We use the Rothe method, i.e.\ the fully implicit time discretisation
  with an equidistant partition of the time interval $I$ with the
  time step $\tau>0$. We denote by $\vvk$, $\Eek$, ... the approximate
  values of $\vv$, $\Ee$, ... at time $k\tau$ with $k=1,2,...,T/\tau$.   
We use the notation for the discretised convective time derivative
$$
%(\cdot)\!\mtd{^{}}
\mtdk{(\cdot)}:=\frac{(\cdot)_\tau^k-(\cdot)_\tau^{k-1}}\tau+\big(\vvk\cdot\nabla\big)(\cdot)_\tau^k\,,
$$
i.e.\ e.g.\ $\mtdk\vv$ will mean
$\frac{\vvk{-}\vvkk}\tau+(\vvk\cdot\nabla)\vvk$ etc.
With this notation, we consider the scheme
\begin{subequations}\label{eq-disc}\begin{align}\label{eq:1a-disc}
    &\density
    %\frac{\vvk{-}\vvkk}\tau+(\vvk\cdot\nabla)\vvk
    \mtdk\vv
    =\operatorname{div}\big({\bm S}_\tau^k
    %\partial_{\Ee}^{}\varphi(\Eek,\alphak,\zzk)
    %{\boldsymbol S}
    +{\boldsymbol S}_{\rm str,\tau}^k
    \!+
    %-\operatorname{div}%{\boldsymbol S}_{\rm hpr,\tau}^k
    k_{\rm v}\EE(\vvk)\big)
    +{\bm f}_\tau^k-
    \frac{\density}2(\operatorname{div}\vvk)\vvk\,,
    \\[.2em]\label{eq:1b-disc}
  &%\bdot
%\frac{\Eek{-}\Eekk}\tau+(\vvk\cdot\nabla)\Eek
\mtdkEe
    =\EE(\vvk)-
    %\frac{\Epk{-}\Epkk}\tau-(\vvk\cdot\nabla)\Epk
    \mtdkEp+k_{\rm e}\Delta{\bm S}_\tau^k%\Eek
  \ \ \ \text{ with }\ \ {\bm S}_\tau^k=\partial_{\Ee}^{}\varphi(\Eek,\alphak,\zzk)
\,,
  %\\
      %  &{\bm S}=\varphi'_{\Ee}(\Ee,\alpha)+\Big(\frac \varrho 2 |\bm v|^2+\varphi(\Ee,\alpha,\zz)+\frac {{\color{red}\text{\huge$\kappa$}}_1}2|\nabla\Ee|^2+\frac {{\color{red}\text{\huge$\kappa$}}_2}2|\nabla\Ep|^2+\frac {{\color{red}\text{\huge$\kappa$}}_3}2|\nabla\alpha|^2\Big){\bm I}+\frac {{\color{red}\text{\huge$\kappa$}}_1}2 \nabla\Ee\otimes\nabla..
 \\\label{eq:3-disc}
 &\partial_{\mtd\Ep}\zeta\Big(\alphakk,\zzkk;
 %\frac{\Epk{-}\Epkk}\tau+(\vvk\cdot\nabla)\Epk
 \mtdkEp
 ,%\mtd\alpha
 %\frac{\alphak{-}\alphakk}\tau+(\vvk\cdot\nabla)\alphak
 \mtdk\alpha\Big)
 -{\bm S}_\tau^k
 %\partial_{\Ee}^{}\varphi(\Eek,\alphak,\zzk)
 %{\boldsymbol S}_{_{\rm KV}}
 \ni%k_{\rm v}\Delta\mtdEp
 k_{\rm p}\Delta\Epk\,,
        \\[-.3em]\label{eq:4-disc}
        &\partial_{\mtd\alpha}\zeta\Big(\alphakk,\zzkk;
        %\mtdEp,
        %\mtd\alpha
        %\frac{\Epk{-}\Epkk}\tau+(\vvk\cdot\nabla)\Epk,
        \mtdkEp,
        %\frac{\alphak{-}\alphakk}\tau+(\vvk\cdot\nabla)\alphak
        \mtdk\alpha\Big)
        +\frac{\alphak{-}\alphakk}{\sqrt\tau}+\partial_\alpha^{}\varphi(\Eek,\alphak,\zzk)\ni k_{\rm a}\Delta\alphak\,,
 \\[.1em]\label{eq:5-disc}
 &%\frac{\zzk{-}\zzkk}\tau+(\vvk\cdot\nabla)\zzk
\mtdk{\zz}
 =\operatorname{div}\big(\mathbb M(\alphakk,\zzkk)\nabla\muk\big)
 %\,{\bm j}(\alpha,\zz,\nabla\mu)
 \ \ \ \text{ with }\ \ \
%\\[.1em]\label{eq:6-disc}&
 \muk=\partial_{\zz}^{}\varphi(\Eek,\alphak,\zzk)\,,
\end{align}
and with the discrete structural stress
\begin{align}%\nonumber
 &{\boldsymbol S}_{\rm str,\tau}^k=k_{\rm p}\nabla\Epk\boxtimes\nabla\Epk+
     k_{\rm a}\nabla\alphak\boxtimes\nabla\alphak
%\\[-.2em]&\qquad\quad\
     -\Big(\varphi(\Eek,\alphak,\zzk)+
  %\frac {k_{\rm e}} 2 |\nabla\Ee|^2+
  \frac{k_{\rm p}}2 |\nabla\Epk|^2
  +\frac{k_{\rm a}}2 |\nabla\alphak|^2
 % +\frac{{\color{red}\text{\huge$\kappa$}}}2 |\nabla\zzk|^2
                                        %-\frac\varrho2|{\bm v}|^2
   \Big){\boldsymbol I}\,.\label{structural-stress-k}
  \end{align}\end{subequations}
%Let us emphasize that the structural stress 
%\eqref{structural-stress-k} couples the whole system \eqref{eq-disc},
%which otherwise would decouple 
%$(\vvk,\Eek,\Epk)$ from $\alphak$ and from $(\zzk,\muk)$.
The boundary conditions \eqref{BC} are discretised correspondingly,
i.e.%\COMMENT{I MODIFIED \eqref{BC2-disc} TO CORRESPOND \eqref{BC2}, OK?}
\begin{subequations}\label{BC-disc}
\begin{align}\label{BC1-disc}
  &
  \vvk\cdot{\bm n}=0,
  \\\label{BC2-disc}
  &
  %\mathbf P_{\textsc t}
  \TTT\big(\EEE(%\partial_{\Ee}^{}\varphi(\Eek,\alphak,\zzk)
  {\bm S}_\tau^k{+}k_{\rm v}\bm e(\vvk){+}{\boldsymbol S}_{\rm str,\tau}^k 
      %{-}\operatorname{div}(k_{\rm v}\nabla\bm e(\vvk)
      ){\bm n}\TTT\big)_{\mathsf t}\EEE
  %-\divS(%{\boldsymbol S}_{\rm hpr}
%k_{\rm v}\nabla\EE(\vvk)
      %{\bm n})
      +\gamma(\vvk)_{\mathsf t}=(\bm g_\tau^k)_{\mathsf t},
  \\
  \label{BC3-disc}
  &
 % \mathbf P_{\textsc t}(
  %{\boldsymbol S}_{\rm hpr}
  %\nabla\EE(\vvk):({\bm n}\otimes{\bm n})
  %%)
  %=\bm 0, \ \ \ \ \ \ \
  ({\bm n}{\cdot}\nabla)\bm S_\tau^k=0,\ \ \ \ \ \ \
  ({\bm n}{\cdot}\nabla)\Epk=0,\ \ \ \ \ \ \ 
 % \\&
  \nabla\alphak{\cdot}{\bm n}=0,\\
  &
  \mathbb M(\alphakk,\zzkk):(\nabla\muk\otimes{\bm n})
  %        {\bm j}(\alpha,\zz,\nabla\mu)\cdot{\bm n}
          =h_\tau^k%\mu\varkappa\mu_{\rm ext,\tau}^k
  \qquad\text{ and }\qquad\nabla\zzk\cdot{{\bm n}}=0\,.\label{BCe-disc}
\end{align}
\end{subequations}
We used the notation $\bm f_\tau^k:=\int_{(k-1)\tau}^{k\tau}\bm f(t)\,\d t$ and
similarly also for $\bm g_\tau^k$ and $h_\tau^k$.
The system of boundary-value problems \eqref{eq-disc}--\eqref{BC-disc}
is to be solved recursively for $k=1,2,...,T/\tau$, assuming
$T/\tau$ integer, and starting with
\begin{align}\label{IC-disc}
\vv_\tau^0=\vv_0,\ \ \ \ \ {\Ee}_{\tau}^0={\Ee}_0,\ \ \ \ \ 
{\Ep}_{\tau}^0={\Ep}_0,\ \ \ \ \ 
\alpha_\tau^0=\alpha_0,\ \ \ \ \ \zz_\tau^0=\zz_0\,.
\end{align}
Let us point out that the term $(\alphak{-}\alphakk)/\sqrt\tau$ in
\eqref{eq:4-disc} is devised to convexify $\varphi$ \SSS using \EEE %due to
\eqref{semi-convex} for small $\tau>0$ %and
\SSS but it still \EEE vanishes in the limit.

\SSS For a given $(\vvkk,\Eekk,\Epkk,\alphakk,\zzkk
%,\muk
)\in H^1(\Omega;\R^d{\times}\R^{d\times d}{\times}\R^{d\times d}{\times}\R^\ell{\times}\R)=: V$, \EEE
the existence of weak solutions $(\vvk,\Eek,\Epk,\alphak,\zzk
%,\muk
)\in V\EEE$ of the coupled semi-linear boundary-value problem
\eqref{eq-disc}--\eqref{BC-disc} can thus be seen by the
\SSS application of Galerkin-approximation-based arguments from the \EEE
theory of coercive
%pseudomonotone
\SSS weakly continuous \EEE set-valued operators \SSS from $V$ to
$Z^*$ for some $Z\subset V$ with the set-valued part arising from a
convex potential; cf.\ e.g.\ \cite[Sect.\,2.5 and 5.3]{Roub13NPDE}.
Here one should choose $Z=W^{1,\infty}(\Omega;\R^d)\times
H^1(\Omega;\R^{d\times d}{\times}\R^{d\times d}{\times}\R^\ell{\times}\R)$ to handle
the structural stress which belongs to $L^1(\Omega;\R^{d\times d})\subset
W^{1,\infty}(\Omega;\R^{d\times d})^*$ but not to
$H^1(\Omega;\R^{d\times d})^*$ in general.
\EEE
 \FFF (Note that the system does not have any potential
because of the convective terms and the related structural stress occuring in
\eqref{eq:1a-disc} make the system nonsymmetric, so that
the direct method cannot be used.) \EEE
The mentioned coercivity is a particular consequence of the \SSS a priori \EEE
estimates derived below. \SSS The weak continuity actually makes the components
$\Epk$ and $\alphak$ strongly convergent by the arguments like
\eqref{strong-nabla} below, which is needed for the continuity of the
nonlinearly dependent structural stress. Also the classical Relich
compact-embedding theorem is used at several places to coup
with the lower-order nonlinearities \EEE
The $L^2$-information about gradients of $\Eek$ and $\zzk$ can be
obtained like in \eqref{est:2++} below.
Moreover, we can rely also on an $L^2$-information about 
$\Delta\Epk$ and $\Delta\alphak$ like in \eqref{est-Delta} below.
Thus the equation/inclusions (\ref{eq-disc}b,c,d) hold even pointwise
a.e.\ on $\Omega$. % Even, we have here also
\SSS We have here additionally \EEE
$\nabla\bm S_\tau^k\in L^2(\Omega;\R^{d\times d})$. 

\medskip\noindent{\it Step 2. (Energetics of the discrete solutions)}:
The \SSS a priori \EEE estimation is based on the energy test. This means here 
the test of \eqref{eq:1a-disc} by $\vvk$
while using also \eqref{eq:1b-disc}, then % the
\SSS we \EEE test the inclusion
\eqref{eq:3-disc} by $(\Epk{-}\Epkk)/\tau+(\vvk{\cdot}\nabla)\Epk$ and the
inclusion
\eqref{eq:4-disc} by $(\alphak{-}\alphakk)/\tau+(\vvk{\cdot}\nabla)\alphak$,
and % the tests of
\SSS we test \EEE
the particular equations in
\eqref{eq:5-disc} %and \eqref{eq:6-disc}
by $\muk$ and $(\zzk{-}\zzkk)/\tau+\vvk{\cdot}\nabla\zzk$\SSS, respectively\EEE.

The mentioned tests thus give the energy balance
\eqref{energy} written as an inequality for the time-discrete approximation.
More specifically, the terms related to inertia in \eqref{eq:1a-disc} uses
the calculus %leading to this energetics then uses 
\begin{align}\nonumber
  %(\varrho\pdt\vvk-f)\cdot\vvk&=
  \Big(\density\frac{\vvk{-}\vvkk}\tau%\frac{\pl\vvk}{\pl t}
  +\density(\vvk{\cdot}\nabla)\vvk-{\bm f}_{{\rm str},\tau}^k\Big){\cdot}\vvk&=
  %\varrho\frac{\pl\vvk}{\pl t}\cdot\vvk
  %\pdt{}\Big(\frac\varrho2|\vvk|^2\Big)
\frac{\density}2\frac{|\vvk|^2-|\vvkk|^2}\tau
  +\density(\vvk{\cdot}\nabla)\vvk\cdot\vvk%+\varrho(\nabla\vvk)^\top\!\vvk\cdot\vvk
  \\&\qquad+\frac{\density}2({\rm div}\,\vvk)|\vvk|^2
  +\tau\frac{\density}2\Big|\frac{\vvk{-}\vvkk}\tau\Big|^2
%\\&\!\stackrel{\eqref{difference}}{=}
%\pdt{}\Big(\frac\varrho2|\vvk|^2\Big)+\varrho(\vvk{\cdot}\nabla)\vvk{\cdot}\vvk+
     %%{\rm div}(\varrho\vvk\otimes\vvk)\cdot\vvk
%\frac\varrho2\vvk{\cdot}%({\rm div}(|v|^2\bbI)
%\nabla|\vvk|^2+\varrho({\rm div}\,\vvk)|\vvk|^2,
\label{test-of-convective}\end{align}
with the ``structural'' force ${\bm f}_{{\rm str},\tau}^k:=-
\frac12\density(\operatorname{div}\vvk)\vvk$, cf.\ the last term in
\eqref{eq:1a-disc}.
%\eqref{structural-force}.
This holds pointwise, and, when integrated over $\Omega$, we further use
%Green's formula $\int_\Omega\frac\varrho2\vvk\cdot\nabla|\vvk|^2\,\d x
%=-\int_\Omega\frac\varrho2({\rm div}\,\vvk)|\vvk|^2\,\d x$ and also
%\eqref{convective-tested}
also \begin{align}%\nonumber
 \int_\Omega\density(\vvk\cdot\nabla)\vvk\cdot \vvk\,\d x&=
  -\int_\Omega\frac{\density}2|\vvk|^2({\rm div}\,\vvk)\,\d x
  +\int_\Gamma\frac{\density}2|\vvk|^2(\vvk\cdot{\bm n})\,\d S
%\\&=-\int_\Omega\Big(\frac{\density}2|\vvk|^2{\bm I}\Big):\EE(\vvk)\,\d x +\int_\Gamma\frac{\density}2|\vvk|^2(\vvk\cdot{\bm n})\,\d S
  \,.
\label{convective-tested}\end{align}
The last term in \eqref{test-of-convective} is non-negative and will simply be
forgotten, which %will give
gives \SSS a discrete analog of \eqref{eq:12} as \EEE the inequality
\begin{align}\nonumber
  %(\varrho\pdt\vvk-f)\cdot\vvk&=
  &\int_\Omega\!\Big(\density\frac{\vvk{-}\vvkk}\tau%\frac{\pl\vvk}{\pl t}
  +\density(\vvk{\cdot}\nabla)\vvk-{\bm f}_{{\rm str},k}\Big)\cdot\vvk\,\d x
  \stackrel{\eqref{test-of-convective}}{\ge}\int_\Omega\!\Big(
\frac{\density}2\frac{|\vvk|^2-|\vvkk|^2}\tau
  %\frac{\d}{\d t}\int_\Omega\frac\varrho2|\vvk|^2\,\d x
  \\&\quad+
  \density(\vvk{\cdot}\nabla)\vvk\cdot\vvk+\frac\varrho2({\rm div}\,\vvk)|\vvk|^2\Big)\,\d x
  \stackrel{\eqref{convective-tested}}{=}\int_\Omega\frac{\density}2\frac{|\vvk|^2-|\vvkk|^2}\tau\,\d x+\int_\Gamma\frac\varrho2|\vvk|^2(\vvk{\cdot}\bm n)\,\d S\,.
\label{test-of-convective+}\end{align}
The last term vanishes due to the boundary condition \eqref{BC1-disc}.
The further term in \eqref{eq:1a-disc} uses the calculus 
\begin{align}\nonumber
  &\int_\Omega\operatorname{div}{\bm S}_\tau^k
  %\partial_{\Ee}^{}\varphi(\Eek,\alphak,\zzk)
  \cdot\vvk\,\d x=\int_\Gamma{\bm S}_\tau^k
  %\partial_{\Ee}^{}\varphi(\Eek,\alphak,\zzk)
  :(\vvk\otimes\bm n)\,\d S
  %\\[-.3em]&\nonumber\hspace*{22em}-
  -\int_\Omega{\bm S}_\tau^k
  %\partial_{\Ee}^{}\varphi(\Eek,\alphak,\zzk)
  :\EE(\vvk)\,\d x
\\&\nonumber\stackrel{\eqref{eq:1b-disc}}{=}
\int_\Gamma%\partial_{\Ee}^{}\varphi(\Eek,\alphak,\zzk)
    {\bm S}_\tau^k:(\vvk\otimes\bm n)\,\d S
%\\[-.3em]&\nonumber\hspace*{3em}-
-    \int_\Omega%\partial_{\Ee}^{}\varphi(\Eek,\alphak,\zzk)
        {\bm S}_\tau^k:\mtdk\Ee
        -%\partial_{\Ee}^{}\varphi(\Eek,\alphak,\zzk)
        {\bm S}_\tau^k:\mtdk\Ep \SSS - \EEE k_{\rm e}{\bm S}_\tau^k
%\partial_{\Ee}^{}\varphi(\Eek,\alphak,\zzk)
:\Delta{\bm S}_\tau^k
%\partial_{\Ee}^{}\varphi(\Eek,\alphak,\zzk)
\,\d x
\\&\nonumber\ \ \ =\int_\Gamma{\bm S}_\tau^k
%\partial_{\Ee}^{}\varphi(\Eek,\alphak,\zzk)
:(\vvk\otimes\bm n)
\SSS + \EEE k_{\rm e}%\partial_{\Ee}^{}\varphi(\Eek,\alphak,\zzk)
{\bm S}_\tau^k:(\bm n\cdot\nabla){\bm S}_\tau^k
\,\d S
\\[-.3em]&%\nonumber
\hspace*{3em}-
\int_\Omega
%\bigg(
    {\bm S}_\tau^k%\partial_{\Ee}^{}\varphi(\Eek,\alphak,\zzk)
:\frac{\Eek{-}\Eekk}\tau
+%\partial_{\Ee}^{}\varphi(\Eek,\alphak,\zzk)
{\bm S}_\tau^k:(\vvk\cdot\nabla)\Eek
%\\&\hspace*{5em}
-%\partial_{\Ee}^{}\varphi(\Eek,\alphak,\zzk)
{\bm S}_\tau^k:\mtdk\Ep \SSS + \EEE k_{\rm e}|\nabla{\bm S}_\tau^k|^2
%\nabla\partial_{\Ee}^{}\varphi(\Eek,\alphak,\zzk)\Vdots\nabla\Eek
%\bigg)
\,\d x\,,
\label{test-by-v}\end{align}
where we abbreviated ${\bm S}_\tau^k=\partial_{\Ee}^{}\varphi(\Eek,\alphak,\zzk)$.
Finally, we have
\begin{align}
  &\int_\Omega \big(\operatorname{div}(\bm S_{\rm str,\tau}^k
  %-\operatorname{div}(k_{\rm v}\nabla\bm e(\vvk))
  +k_{\rm v}\bm e(\vvk)\big)\cdot\bm v_\tau^k\,\d x\nonumber
\\
&\qquad =\int_\Gamma(\bm S_{\rm str,\tau}^k
+k_{\rm v}\bm e(\vvk)
%-\operatorname{div}(k_{\rm v}\nabla\bm e(\vvk))
):(\vvk\otimes\bm n)\,\d S-\int_\Omega(\bm S_{\rm str,\tau}^k
+k_{\rm v}\bm e(\vvk)
%-k_{\rm v}\Delta\bm e(\vvk)
):\nabla\vvk\,\d x\,.\label{eq:2+}
\end{align}
% Since $\Gamma$ is a smooth surface with no boundary, the application of the divergence theorem on $\Gamma$ yields the last of the above equalities.

The mentioned test \eqref{eq:3-disc}
by $\mtdk{}\Ep=\frac{\Epk-\Epkk}\tau+(\vvk{\cdot}\nabla)\Epk$ gives
%a contribution both
\SSS contributions \EEE
to the dissipation rate and to the stored-energy rate.
The dissipation and the gradient terms in \eqref{eq:3-disc} yield,
\SSS using also a discrete version of
the calculus behind \eqref{eq:14}, that \EEE
\begin{align}\nonumber&\int_\Omega
\partial_{\mtd\Ep}\zeta\Big(\alpha\TTT_\tau^{k-1}\EEE,\zz\TTT_\tau^{k-1}\EEE);\mtdk\Ep,\mtdk\alpha\Big):\mtdk\Ep-
%k_{\rm v}\Delta{\bm{\varPi}}
k_{\rm p}\Delta\Epk:\mtdk\Ep\,\d x
 \\[-.3em]&\nonumber\hspace*{2em}\ge
 \int_\Omega\partial_{\mtd\Ep}\zeta\Big(\alpha\TTT_\tau^{k-1}\EEE,\zz\TTT_\tau^{k-1}\EEE);\mtdk\Ep,\mtdk\alpha\Big)
 :\mtdk\Ep \SSS - \EEE k_{\rm p}(\vvk\cdot\nabla)\Epk:\Delta\Epk\,\d x
 \\[-.3em]&\nonumber\hspace*{7em}+
% \frac{\d}{\d t}
 \int_\Omega\frac{k_{\rm p}}2%|\nabla\Epk|^2
 \frac{|\nabla\Epk|^2-|\nabla\Epkk|^2}\tau\,\d x
-\int_\Gamma k_{\rm p}\nabla\Epk\Vdots\Big(\mtdk\Ep\otimes\bm n\Big)\,\d S
\\[-.3em]&\nonumber\hspace*{2em}=
\int_\Omega
\partial_{\mtd\Ep}\zeta\Big(\alpha\TTT_\tau^{k-1}\EEE,\zz\TTT_\tau^{k-1}\EEE);\mtdk{\Ep\!\!},\mtdk\alpha\Big){:}\mtdk{\Ep\!\!}
+\Big(k_{\rm p}\nabla\Epk\boxtimes\nabla\Epk-
\frac{k_{\rm p}}2|\nabla\Epk|^2\bm I\Big){:}\bm e(\vvk)\,\d x
\\[-.1em]&\hspace*{25em}
+%\frac{\d}{\d t}
\int_\Omega\frac{k_{\rm p}}2\frac{|\nabla\Epk|^2-|\nabla\Epkk|^2}\tau\,\d x\,,
\label{calc-plastic}\end{align}
where the term boundary term
$k_{\rm p}\nabla\Epk\Vdots(\mtdk{\Ep\!\!}\otimes\bm n)=(\bm n\cdot\nabla)\Epk:\mtdk\Ep$ vanishes thanks to \eqref{BC3-disc} and where the inequality relies on the convexity 
of the functional $\Ep\mapsto\int_\Omega \frac{k_{\rm p}}2|\nabla\Ep|^2\,\d x$.
\SSS The inequality follows from the calculus
\begin{align}\nonumber-\int_\Omega\Delta\Epk:\frac{\Epk{-}\Epkk\!\!\!\!}\tau\ \d x
&=\int_\Omega\nabla\Epk:\nabla\frac{\Epk{-}\Epkk\!\!\!\!}\tau\ \d x
\\&\nonumber=\int_\Omega\bigg(\frac{|\nabla\Epk|^2-|\nabla\Epkk|^2\!\!}{2\tau}
\\[-.3em]\nonumber
&\qquad+\frac\tau2\Big|\frac{\nabla\Epk-\nabla\Epkk}{\tau}\Big|^2\bigg)\,\d x
\ge\int_\Omega\frac{|\nabla\Epk|^2-|\nabla\Epkk|^2\!\!}{2\tau}\ \d x.
\end{align}
while the meaning of the expression
$\partial_{\mtd\Ep}\zeta(\alpha,\zz;\DT\Ep,\DT\alpha):\DT\Ep$ is well defined
even if $\partial_{\mtd\Ep}\zeta(\alpha,\zz;\cdot,\DT\alpha)$ is multivalued
at $\DT\Ep=0$, cf.\ \eqref{smoothnes-1}.\EEE
%For the last equality in \eqref{calc-plastic}, we used %a calculation similar to that performed in
%\eqref{eq:14} to transform the integral of the convective term $(\vvk\cdot\nabla)\Epk:\Delta\Epk$.

Moreover, the test of \eqref{eq:4-disc} by
$\mtdk{}\alpha=\frac{\alphak-\alphakk}\tau+\vvk\cdot\nabla\alphak$ gives rise
to the term $\int_\Omega(\vvk\cdot\nabla\alphak)\Delta\alphak\,\d x$. By proceeding as in \eqref{calc-plastic}, using the boundary condition
$\vvk{\cdot}\bm n=0$, we obtain
\begin{align}\nonumber
  &\int_\Omega
\partial_{\mtd\alpha}\zeta\Big(\alpha\TTT_\tau^{k-1}\EEE,\zz\TTT_\tau^{k-1}\EEE);\mtdk\Ep,\mtdk\alpha\Big)\cdot\mtdk\alpha
%k_{\rm v}\Delta{\bm{\varPi}}
+\Big(\partial_\alpha^{}\varphi(\Eek,\alphak,\zzk)+\frac{\alphak{-}\alphakk\!\!\!}{\sqrt\tau}%\Big)\cdot\mtdk\alpha
-k_{\rm a}\Delta\alphak\Big)\cdot\mtdk\alpha\,\d x
 \\[-.3em]&\nonumber\hspace*{2em}\ge
\int_\Omega
\partial_{\mtd\alpha}\zeta\Big(\alpha\TTT_\tau^{k-1}\EEE,\zz\TTT_\tau^{k-1}\EEE);\mtdk\Ep,\mtdk\alpha\Big)\cdot\mtdk\alpha
+\Big(k_{\rm a}\nabla\alphak\boxtimes\nabla\alphak-
\frac{k_{\rm a}}2|\nabla\alphak|^2\bm I\Big):\bm e(\vvk)\,\d x
\\[-.1em]&\hspace*{12em}
+%\frac{\d}{\d t}
\int_\Omega\partial_\alpha^{}\varphi(\Eek,\alphak,\zzk)\cdot\frac{\alpha_\tau^k-\alpha^{k-1}_\tau\!\!\!}\tau+\frac{k_{\rm a}}2\frac{|\nabla\alphak|^2-|\nabla\alphakk|^2}\tau\,\d x\,.
\label{calc-alpha}\end{align}
%\COMMENT{HERE MAYBE TO CHECK}
\SSS The inequality in \eqref{calc-alpha} arises from the same reasons as
in \eqref{calc-plastic} using \eqref{smoothnes-2}.\EEE

Eventually,
%for
the test of \eqref{eq:5-disc} by $\muk$
%being a weak solution to \eqref{eq:6-disc} gives rise to the
%term $\int_\Omega (\vvk\cdot\nabla\zzk)\Delta\zzk\,\d x$. For this term,
%we use again the above calculus with $\zz$ instead of $\alpha$.
%Also, we use \eqref{eq:6-disc} tested by $\frac{\zzk-\zzkk}\tau$, which
gives
\begin{align}\nonumber
  &\int_\Omega\Big(\mtdk\zz
  -{\rm div}(\mathbb M(\alphakk,\zzkk)\nabla\muk)\Big)\,\muk
%\,{\bm j}(\alpha\TTT_\tau^{k-1}\EEE,\zz\TTT_\tau^{k-1}\EEE),\nabla\mu_k)
\,\d x
\\&\nonumber=\int_\Omega\Big(\frac{\zzk{-}\zzkk}\tau
+\vvk\cdot\nabla\zzk\Big)\muk+\mathbb M(\alphakk,\zzkk)\nabla\muk
%{\bm j}(\alpha\TTT_\tau^{k-1}\EEE,\zz\TTT_\tau^{k-1}\EEE),\nabla\mu_k)
{\cdot}\nabla\muk\,\d x-\int_\Gamma h\muk
%{\bm j}(\alpha\TTT_\tau^{k-1}\EEE,\zz\TTT_\tau^{k-1}\EEE),\nabla\mu_k)
\,\d S
%\\&\nonumber=
%\int_\Omega\partial_\zz^{}\varphi(\Eek,\alphak,\zzk)\Big(\frac{\zzk{-}\zzkk}\tau
%+\vvk\cdot\nabla\zzk\Big)+\mathbb M(\alphakk,\zzkk)\nabla\muk
%%{\bm j}(\alpha\TTT_\tau^{k-1}\EEE,\zz\TTT_\tau^{k-1}\EEE),\nabla\mu_k)
%\cdot\nabla\muk\,\d x
\\&\ge\int_\Omega\!\partial_\zz^{}\varphi(\Eek,\alphak,\zzk)
\Big(\frac{\zzk{-}\zzkk\!\!}\tau+\vvk{\cdot}\nabla\zzk\Big)
+\mathbb M(\alphakk,\zzkk)\nabla\muk
%{\bm j}(\alpha\TTT_\tau^{k-1}\EEE,\zz\TTT_\tau^{k-1}\EEE),\nabla\mu_k)
\cdot\nabla\muk\,\d x \SSS -\int_\Gamma h\muk
%{\bm j}(\alpha\TTT_\tau^{k-1}\EEE,\zz\TTT_\tau^{k-1}\EEE),\nabla\mu_k)
\,\d S\EEE\,.
\label{est-chi}\end{align}
%Here we again used the calculus \eqref{Korteweg-calculus} but now for
%$\zz$ instead of $\Ep$.
%

Using the semi-convexity of $\varphi$, we can estimate the sum of the three
terms arising in \eqref{test-by-v}, \eqref{calc-alpha}, and \eqref{est-chi} together with the convexifying term
in \eqref{eq:4-disc}
as
\begin{align}\nonumber
  &%\int_\Omega
  %\partial_\Ee^{}\varphi(\Eek,\alphak,\zzk)
  {\bm S}_\tau^k:\frac{\Eek{-}\Eekk}\tau
%\\&\nonumber\hspace*{8em}
+\Big(\partial_\alpha^{}\varphi(\Eek,\alphak,\zzk)
  +\frac{\alphak{-}\alphakk}{\sqrt\tau}\Big)\cdot\frac{\alphak{-}\alphakk}\tau
  +\partial_\zz^{}\varphi(\Eek,\alphak,\zzk)\frac{\zzk{-}\zzkk}\tau%\,\d x
\\&\nonumber\hspace*{0em}
=
%\int_\Omega
\partial_\Ee^{}\varphi(\Eek,\alphak,\zzk):\frac{\Eek{-}\Eekk}\tau+
  \Big(\partial_\alpha^{}\varphi(\Eek,\alphak,\zzk)
  +\frac{\alphak}{\sqrt\tau}\Big)\cdot\frac{\alphak{-}\alphakk}\tau
\\&\nonumber\hspace*{8em}
+\partial_\zz^{}\varphi(\Eek,\alphak,\zzk)
\frac{\zzk{-}\zzkk}\tau-\frac{\alphakk}{\sqrt\tau}\cdot\frac{\alphak{-}\alphakk}\tau%+(1{-}\sqrt[4]\tau)\Big|\frac{\alphak{-}\alphakk}\tau\Big|^2
%\,\d x
  \\&\nonumber\hspace*{0em}
  \ge
  %\int_\Omega
  \frac{\varphi(\Eek,\alphak,\zzk)
  -\varphi(\Eekk,\alphakk,\zzkk)}\tau+\frac1{2\sqrt\tau}\frac{|\alphak|^2{-}|\alphakk|^2}\tau
-\frac{\alphakk}{\sqrt\tau}\cdot\frac{\alphak{-}\alphakk}\tau%+(1{-}\sqrt[4]\tau)\Big|\frac{\alphak{-}\alphakk}\tau\Big|^2
%\,\d x
  \\&\hspace*{0em}
  =
  %\int_\Omega
  \frac{\varphi(\Eek,\alphak,\zzk)-\varphi(\Eekk,\alphakk,\zzkk)}\tau
-\frac{\sqrt\tau}2\Big|\frac{\alphak{-}\alphakk}\tau\Big|^2
%\,\d x
\,,\label{semi-convex-trick}
\end{align}
cf.\ also the calculation in \cite[Remark~8.24]{Roub13NPDE}. This holds
a.e.\ on $\Omega$ and is to be integrated over $\Omega$.
For the remaining three convective terms arising from these tests, we use
the calculus
\begin{align}\nonumber&
\int_\Omega\Big(\partial_{\Ee}^{}\varphi(\Eek,\alphak,\zzk):(\vvk\cdot\nabla)\Eek
+\partial_\alpha^{}\varphi(\Eek,\alphak,\zzk)
\cdot(\vvk\cdot\nabla\alphak)
\\[-.3em]\nonumber
&\qquad\qquad
+\partial_\zz^{}\varphi(\Eek,\alphak,\zzk)\cdot(\vvk\cdot\nabla\zzk)\Big)
\,\d x
=\int_\Omega
\nabla\varphi(\Eek,\alphak,\zzk)\cdot\vvk\,\d x
\\[-.3em]&\nonumber\qquad\qquad\qquad\qquad
=\int_\Gamma\varphi(\Eek,\alphak,\zzk)\vvk\cdot{\bm n}\,\d S-\int_\Omega
\varphi(\Eek,\alphak,\zzk){\rm div}\,\vvk\,\d x
\\[-.3em]&\qquad\qquad\qquad\qquad\qquad\qquad
=-\int_\Gamma\varphi(\Eek,\alphak,\zzk)\bm I:\EE(\vvk)\,\d x\,,
\label{calc-energy-pressure}\end{align}
which cancels with the pressure-type stress contribution
$\varphi(\Eek,\alphak,\zzk)\bm I$.
%\COMMENT{CAREFULL -- \eqref{eq+} TO BE USED!!!}

Eventually,  after %integration over a time interval $[0,t]$,
summation over $k=1,2,....$, 
we obtain \eqref{energy} as an upper estimate up to
an error term which is small for $\tau>0$ small, so that it
can be used for \SSS a priori \EEE estimates. More precisely,
by the test of the regularizing term $(\alphak{-}\alphakk)/\sqrt\tau$ by
$\mtdk{}\alpha$, we obtain still the term
\begin{align}\nonumber
\int_\Omega\frac{\alphak{-}\alphakk}{\sqrt\tau}\cdot(\vvk\cdot\nabla)\alphak\,\d x
&=\sqrt\tau\int_\Omega\frac{\alphak{-}\alphakk}{\tau}\cdot(\vvk\cdot\nabla)\alphak\,\d x
\\\nonumber&
=\sqrt\tau\int_\Omega\frac{\alphak{-}\alphakk}{\tau}\cdot\mtdk\alpha
-\Big|\frac{\alphak{-}\alphakk}{\tau}\Big|^2\,\d x
\\&\le
\frac{\sqrt\tau}2\Big\|\mtdk\alpha\Big\|_{L^2(\Omega;\R^\ell)}^2-\frac{\sqrt\tau}2\Big\|\frac{\alphak{-}\alphakk}{\tau}\Big\|_{L^2(\Omega;\R^\ell)}^2\,,
\label{semi-convex-trick+}\end{align}
which allows for estimation in the next step when relying on the
assumption \eqref{ass:2} and on the last term in \eqref{semi-convex-trick}.
%written for the discrete solution

\medskip\noindent{\it Step 3. (\SSS A priori \EEE estimates)}:
Using the values $(\vvk)_{k=0}^{T/\tau}$, we define the piecewise constant and
the piecewise affine interpolants respectively as
%\begin{subequations}\label{def-of-interpolants}
\begin{align}\label{def-of-interpolants}
&\overline\vv_\tau(t):=\vvk,\ \ \ \underline\vv_\tau(t):=\vvkk,
%\\&
\ \ \text{ and }\ \ \vv_\tau(t):=\Big(\frac t\tau{-}k{+}1\Big)\vvk
+\Big(k{-}\frac t\tau\Big)\vvkk\ \ \text{ for }\ \
(k{-}1)\tau<t\le k\tau
\end{align}%\end{subequations}
for $k=0,1,...,T/\tau$. Analogously, we define also 
$\Eetau$, $\overlineEetau$, etc. In terms of such interpolants,
%using also the notation
%$$\mtdtau{(\cdot)}:=\Big[\pdt{}+\overline\vv_\tau\cdot\nabla\Big](\cdot)\,,$$
we can
write the discrete recursive system \eqref{eq-disc} ``compactly'' as
\begin{subequations}\label{eq-disc+}\begin{align}\label{eq:1a-disc+}
    &\density
    %\frac{\vvk{-}\vvkk}\tau+(\vvk\cdot\nabla)\vvk
    %\mtdtau{\vv_\tau}
    \pdt{\vv_\tau}+(\overline\vv_\tau{\cdot}\nabla)\overline\vv_\tau
    =\operatorname{div}\big(\overlineStau
%\partial_{\Ee}^{}\varphi(\overlineEetau,\overline\alpha_\tau,\overline\zz_\tau)
    +\overline{{\boldsymbol S}}_{\rm str,\tau}
+k_{\rm v}\EE(\overline\vv_\tau)
    %-\operatorname{div}%{\boldsymbol S}_{\rm hpr,\tau}^k
    %k_{\rm v}\nabla\EE(\overline\vv_\tau)
    \big)
    +\overline{\bm f}_\tau-
    \frac{\density}2(\operatorname{div}\overline\vv_\tau)\overline\vv_\tau\,,
    \\[.2em]\label{eq:1b-disc+}
    &%\mtdtau{\Eetau}
    \pdt{\Eetau}+(\overline\vv_\tau{\cdot}\nabla)\overlineEetau
    =\EE(\overline\vv_\tau)-
    %\mtdkEpk
    \pdt{\Eptau}-(\overline\vv_\tau{\cdot}\nabla)\overlineEptau
    +k_{\rm e}\Delta\overlineStau
    %\overlineEetau
\ \ \ \text{ with }\ \ \overlineStau=
 \partial_{\Ee}^{}\varphi(\overlineEetau,\overline\alpha_\tau,\overline\zz_\tau)\,,
 \\\label{eq:3-disc+}
 &\partial_{\mtd\Ep}\zeta\Big(\underline\alpha_\tau,\underline\zz_\tau;
 %\frac{\Epk{-}\Epkk}\tau+(\vvk\cdot\nabla)\Epk
 %\mtdkEpk,%\mtdk\overline\alpha_\tau
 \pdt{\Eptau}{+}(\overline\vv_\tau{\cdot}\nabla)\overlineEptau,
 \pdt{\alpha_\tau}{+}(\overline\vv_\tau{\cdot}\nabla)\overline\alpha_\tau
 \Big)
 -\overlineStau
%\partial_{\Ee}^{}\varphi(\overlineEetau,\overline\alpha_\tau,\overline\zz_\tau)
 \ni k_{\rm p}\Delta\overlineEptau\,,
        \\[.2em]\label{eq:4-disc+}
        &\partial_{\mtd\alpha}\zeta\Big(\underline\alpha_\tau,\underline\zz_\tau;
        %\mtdkEpk,
        \pdt{\Eptau}{+}(\overline\vv_\tau{\cdot}\nabla)\overlineEptau,
        %\mtdk\overline\alpha_\tau
\pdt{\alpha_\tau}{+}(\overline\vv_\tau{\cdot}\nabla)\overline\alpha_\tau
        \Big)+\sqrt\tau\pdt{\alpha_\tau}+\partial_\alpha^{}\varphi(\overlineEetau,\overline\alpha_\tau,\overline\zz_\tau)\ni k_{\rm a}\Delta\overline\alpha_\tau\,,
 \\[.1em]\label{eq:5-disc+}
 &%\mtdk{\zz_\tau}
\pdt{\zz_\tau}+(\overline\vv_\tau{\cdot}\nabla)\overline\zz_\tau
 =\operatorname{div}(\mathbb M(\underline\alpha_\tau,\underline\zz_\tau)\nabla\overline\mu_\tau) \ \ \ \text{ with }\ \ \
%\\[.1em]\label{eq:6-disc}&
 \overline\mu_\tau=\partial_{\zz}^{}\varphi(\overlineEetau,\overline\alpha_\tau,\overline\zz_\tau)\,,
\end{align}
and with the discrete structural stress
\begin{align}%\nonumber
  &\overline{{\boldsymbol S}}_{\rm str,\tau}=
  k_{\rm p}\nabla\overlineEptau\boxtimes\nabla\overlineEptau+
     k_{\rm a}\nabla\overline\alpha_\tau\boxtimes\nabla\overline\alpha_\tau
%\\&\qquad
-\Big(\varphi(\overlineEetau,\overline\alpha_\tau,\overline\zz_\tau)+
  %\frac {k_{\rm e}} 2 |\nabla\Ee|^2+
  \frac{k_{\rm p}}2 |\nabla\Epk|^2
  +\frac{k_{\rm a}}2 |\nabla\overline\alpha_\tau|^2
 % +\frac{{\color{red}\text{\huge$\kappa$}}}2 |\nabla\overline\zz_\tau|^2
                                        %-\frac\varrho2|{\bm v}|^2
   \Big){\boldsymbol I}\label{structural-stress+}
  \end{align}\end{subequations}
and with the boundary conditions \eqref{BC-disc} written analogously.
Actually, like \eqref{VI}, the inclusions (\ref{eq-disc+}c,d) mean
\begin{align}\nonumber
&\int_0^T\!\!\!\int_\Omega
  \zeta(\underline\alpha_\tau,\underline\zz_\tau;\widetildeEp,\widetilde{\alpha})-\partial_{\Ee}^{}\varphi(\overlineEetau,\overline\alpha_\tau,\overline\zz_\tau)
  :\Big(\widetildeEp{-}\pdt{\Eptau}{-}(\overline\vv_\tau{\cdot}\nabla)\overlineEptau\Big)
\\[-.3em]\nonumber&\qquad
  +\Big(\partial_\alpha^{}\varphi(\overlineEetau,\overline\alpha_\tau,\overline\zz_\tau){+}\sqrt\tau\pdt{\alpha_\tau}\Big)
      \cdot\Big(\widetilde{\alpha}{-}\pdt{\alpha_\tau}{-}(\overline\vv_\tau{\cdot}\nabla)\overline\alpha_\tau\Big)
 %+k_{\rm p}\nabla\Ep\Vdots\nabla\Big(\widetildeEp{-}\pdt{\Eptau}-(\overline\vv_\tau{\cdot}\nabla)\overlineEptau\Big)
    + k_{\rm p}\nabla\overlineEptau\Vdots\nabla\widetildeEp \\[-.3em]\nonumber&\qquad
+ k_{\rm p}\Delta\overlineEptau{:}(\vv{\cdot}\nabla)\overlineEptau
+k_{\rm a}\nabla\overline\alpha_\tau{:}\nabla\widetilde{\alpha}
+k_{\rm a}\Delta\overline\alpha_\tau{\cdot}(\vv{\cdot}\nabla)\overline\alpha_\tau
      %+k_{\rm a}\nabla\alpha:\nabla\Big(\widetilde{\alpha}{-}\pdt{\alpha_\tau}{-}(\overline\vv_\tau{\cdot}\nabla)\overline\alpha_\tau\Big)
      \,\d x\d t
+\int_\Omega\frac{k_{\rm p}}2|\nabla\Ep_0|^2+\frac{k_{\rm p}}2|\nabla\alpha_0|^2\,\d x
 \\[-.3em]&\nonumber
 \ge\int_\Omega\!\frac{k_{\rm p}}2|\nabla\Ep_\tau(T)|^2+\frac{k_{\rm p}}2|\nabla\alpha_\tau(T)|^2\,\d x
 \\[-.3em]&\qquad+\int_0^T\!\!\!\int_\Omega\zeta\Big(\underline\alpha_\tau,\underline\zz_\tau;\pdt{\Eptau}{+}(\overline\vv_\tau{\cdot}\nabla)\overlineEptau,\pdt{\alpha_\tau}{+}(\overline\vv_\tau{\cdot}\nabla)\overline\alpha_\tau\Big)\,\d x\d t
   \label{VI-disc}
\end{align}
for any
$(\widetildeEp,\widetilde{\alpha})\in L^2(I;H^1(\Omega;\R^{d\times d}{\times}\R^{\ell}))$.
%\COMMENT{PASAGE FROM $\zeta$ TO $\xi$ IS NOW NEW:}
%Let us note that, from the inclusion $\partial%_{(\mtd{\Ep},\mtd\alpha)}
% \zeta_{\alpha,\zz}^{(1)}(\mtd{\Ep})\ni\bm{F}$, i.e.\
% from the inequality $\zeta_{\alpha,\zz}^{(1)}(\wtEp)
%\ge \bm{F}{:}(\wtEp{-}\mtd{\Ep})+\zeta_{\alpha,\zz}^{(1)}(\mtd{\Ep})$,
%by putting $\wtEp=0$, we can only obtain the
% inequality $\zeta_{\alpha,\zz}^{(1)}(\mtd{\Ep})\le \bm{F}{:}\mtd{\Ep}$,
%  which is not exactly the expected energy (im)balance, although it
%  would suffice for derinving the \SSS a priori \EEE estimates.
  %For proving
  %the energetics \eqref{energy} later in Step~5, let us derive the expected
  %discrete analog of \eqref{energy} even now. As $\zeta_{\alpha,\zz}^{(1)}$
  %is differentiable except 0, cf.\ \eqref{ass:2}, from
  %$\partial%_{(\mtd{\Ep},\mtd\alpha)}
  %\zeta_{\alpha,\zz}^{(1)}(\mtd{\Ep})\ni\bm{F}$, i.e.\
  %from $\zeta_{\alpha,\zz}^{(1)}(\wtEp)
  %\ge \bm{F}{:}\wtEp{-}\mtd{\Ep})+
  %\zeta_{\alpha,\zz}^{(1)}(\mtd{\Ep})$, by using a test
  %$\wtEp=\mtd{\Ep}\pm\epsilon\whEp$
  %with $\epsilon\to0$ (if $\mtd{\Ep}\ne0$) and by puting eventually 
%$\whEp=\mtd{\Ep}$,
%    we can
%    conclude that $\partial\zeta_{\alpha,\zz}^{(1)}(\mtd{\Ep}){:}\mtd{\Ep}
%    =\bm{F}{:}\mtd{\Ep}$. This identity obviously holds also for
%    $\mtd{\Ep}=0$, although the argumentation does not.
%    The same arguments apply for $\zeta_{\alpha,\zz}^{(2)}$ from the decomposition
%    \eqref{ass:2}.
 % Thus, from (\ref{eq-disc+}c,d) which actually holds
 %  pointwise a.e.\ on $I{\times}\Omega$ because
Let us note that
%we have also
$\Delta\overlineEptau\in L^2(I{\times}\Omega;\R^{d\times d})$
  and $\Delta\overline\alpha_\tau\in L^2(I{\times}\Omega;\R^\ell)$,
  so that the integrals in \eqref{VI-disc} have a good sense.
By putting $\widetildeEp=0$ and $\widetilde{\alpha}=0$, from \eqref{VI-disc}
we can also read
\begin{align}\nonumber
&\int_0^T\!\!\!\int_\Omega
\partial_{\Ee}^{}\varphi(\overlineEetau,\overline\alpha_\tau,\overline\zz_\tau)
  :\Big(\pdt{\Eptau}{+}(\overline\vv_\tau{\cdot}\nabla)\overlineEptau\Big)
%\\[-.3em]\nonumber&\qquad
  -\Big(\partial_\alpha^{}\varphi(\overlineEetau,\overline\alpha_\tau,\overline\zz_\tau){+}\sqrt\tau\pdt{\alpha_\tau}\Big)
      \cdot\Big(\pdt{\alpha_\tau}{+}(\overline\vv_\tau{\cdot}\nabla)\overline\alpha_\tau\Big)
      \\[-.5em]\nonumber&\qquad\qquad\qquad
+ k_{\rm p}\Delta\overlineEptau{:}(\vv{\cdot}\nabla)\overlineEptau
+k_{\rm a}\Delta\overline\alpha_\tau{\cdot}(\vv{\cdot}\nabla)\overline\alpha_\tau
      %+k_{\rm a}\nabla\alpha:\nabla\Big(\widetilde{\alpha}{-}\pdt{\alpha_\tau}{-}(\overline\vv_\tau{\cdot}\nabla)\overline\alpha_\tau\Big)
      \,\d x\d t
+\int_\Omega\frac{k_{\rm p}}2|\nabla\Ep_0|^2+\frac{k_{\rm p}}2|\nabla\alpha_0|^2\,\d x
 \\[-.3em]&%\nonumber
 \ge\int_0^T\!\!\!\int_\Omega\zeta\Big(\underline\alpha_\tau,\underline\zz_\tau;\pdt{\Eptau}{+}(\overline\vv_\tau{\cdot}\nabla)\overlineEptau,\pdt{\alpha_\tau}{+}(\overline\vv_\tau{\cdot}\nabla)\overline\alpha_\tau\Big)\,\d x\d t
 +\int_\Omega\!\frac{k_{\rm p}}2|\nabla\Ep_\tau(T)|^2\!+\frac{k_{\rm p}}2|\nabla\alpha_\tau(T)|^2\,\d x
 %\\[-.3em]&\qquad
 \,.
   \label{VI-disc+}
  \end{align}

Of course, we can write the above estimates on $[0,k\tau]$ with
$k=1,...,T/\tau$ instead of $I=[0,T]$. Altogether, we obtain a discrete
%analog of the
energy-like balance 
\begin{align}\nonumber
&\int_\Omega\frac\varrho2|\vv_\tau(t)|^2+
  \varphi(\Ee_\tau(t),\alpha_\tau(t),\zz_\tau(t))
  +\frac{k_{\rm p}}2|\nabla\Ep_\tau(t)|^2+\frac{k_{\rm a}}2|\nabla\alpha_\tau(t)|^2
  \,\d x
%\\&\nonumber\qquad
  \\&\qquad\quad\nonumber
  +\int_0^t\!\!\bigg(\int_\Omega
  \zeta\Big(\underline\alpha_\tau,\underline\zz_\tau;\pdt{\Eptau}{+}(\overline\vv_\tau{\cdot}\nabla)\overlineEptau,\pdt{\alpha_\tau}{+}(\overline\vv_\tau{\cdot}\nabla)\overline\alpha_\tau\Big)  +k_{\rm v}|\EE(\overline\vv_\tau)|^2
 \\&\qquad\qquad\nonumber
  +{\mathbb M}(\underline\alpha_\tau,\underline\zz_\tau)\nabla\overline\mu_\tau{\cdot}\nabla\overline\mu_\tau
+k_{\rm e}\big|\nabla\overlineStau\big|^2+{\sqrt\tau}\Big|\pdt{\alpha_\tau}\Big|^2
  \,\d x
  +\int_\Gamma\gamma|\overline\vv_{\rm t,\tau}|^2
  \\&\qquad\qquad\quad\nonumber
  \le\int_0^t\!\!\bigg(\int_\Omega\overline{\bm f}_\tau\cdot\overline\vv_\tau
+\frac{\sqrt\tau}2\Big|\pdt{\alpha_\tau}{+}(\overline\vv_\tau{\cdot}\nabla)\overline\alpha_\tau\Big|^2
%+\frac{\sqrt\tau}2\Big|\frac{\partial\alpha_\tau}{\partial t}+(\overline\vv\cdot\nabla)\overline\alpha_\tau\Big|^{2}
  \,\d x
  +\int_\Gamma\overline{\bm g}_{\rm t,\tau}\cdot\overline\vv_{\rm t,\tau}
  +h\overline\mu_\tau\d S\bigg)\d t
  \\&\hspace{18.5em}
  +\int_\Omega\frac\varrho2|\vv_0|^2+\varphi({\Ee}_0,\alpha_0,\zz_0)
   +\frac {k_{\rm a}} 2 |\nabla\alpha_0|^2\,\d x
  \label{energy+}\end{align}
for any $t=k\tau$. It should be emphasized that, as
\eqref{energy+} involves the dissipation potential $\zeta$
and not the dissipation rate $\zeta$, it is not a direct
discrete analog of the energy balance \eqref{energy}, but it
is sufficient for the \SSS a priori \EEE estimates. In fact,
refining the argumentation, \eqref{energy+} with $\xi$ could
have been proved, too.

From the energetic inequality \eqref{energy+} by \TTT using the Young inequality
for \EEE estimating $\int_\Omega{\bm f}_\tau^k\cdot\vvk\,\d x\le
\SSS\|{\bm f}_\tau^k\|_{L^2(\Omega;\R^d)}^{}\|\vvk\|_{L^2(\Omega;\R^d)}^{}
\le\|{\bm f}_\tau^k\|_{L^2(\Omega;\R^d)}^{}(1+\|\vvk\|_{L^2(\Omega;\R^d)}^2)\EEE$ 
%
%The last term in \eqref{test-by-v} is to be handled by using the
%calculus \eqref{calculus-for-diffusion-Ee}.
%\begin{align}\nonumber&
%  \epsilon\|\nabla\Eek\|_{L^2(\Omega;\R^{d\times d\times d})}^2 \le\int_\Omega
 %k_{\rm e}\partial_{\Ee\Ee}^2\varphi(\Eek,\alphak,\zzk)\nabla\Eek\Vdots\nabla\Eek\,\d x
%\\&\nonumber=\int_\Omega
%  k_{\rm e}\nabla\partial_{\Ee}^{}\varphi(\Eek,\alphak,\zzk)\Vdots\nabla\Eek
%-k_{\rm e}\nabla\partial_{\Ee}^{}\varphi(\Eek,\alphak,\zzk)\Vdots\nabla\Eek
%-k_{\rm e}\partial_{\Ee\zz}^2\varphi(\Eek,\alphak,\zzk)\nabla\zzk\Vdots\nabla\Eek
%\,\d x\\&\nonumber\le\int_\Omega
%k_{\rm e}\nabla\partial_{\Ee}^{}\varphi(\Eek,\alphak,\zzk)\Vdots\nabla\Eek\,\d x
% +\frac{k_{\rm e}}\epsilon\|\partial_{\Ee\alpha}^2\varphi(\Eek,\alphak,\zzk)\nabla\alphak\|_{L^2(\Omega;\R^{d\times d\times d})}^2
% \\[-.3em]&\nonumber\qquad
% +\frac{k_{\rm e}}\epsilon
%\|\partial_{\Ee\zz}^2\varphi(\Eek,\alphak,\zzk)\nabla\zzk\|_{L^2(\Omega;\R^{d\times d\times d})}^2
% -\frac\epsilon2\|\nabla\Eek\|_{L^2(\Omega;\R^{d\times d\times d})}^2
%\end{align}where $\epsilon>0$ is from \eqref{ass:1}.
and by using the discrete Gronwall inequality, we obtain the following \SSS a priori \EEE estimates:
\begin{subequations}\label{est}\begin{align}\label{est:1}
    &\|\overline\vv_\tau\|_{L^\infty(I;L^2(\Omega;\R^d))\,\cap\,L^2(I;H^1(\Omega;\R^d))}^{}\le C,
    \\&\label{est:2}
    %\|\Eek\|_{L^2(I;H^1(\Omega;\R^{d\times d}))}^{}\le C,
    \big\|\overlineEetau\big\|_{L^\infty(I;L^2(\Omega;\R^{d\times d}))
    %  \,\cap\,L^2(I;H^1(\Omega;\R^{d\times d}))
    }^{}\le C
    \ \ \ \text{ and }\ \ \
%\|\mtdEe\|_{L^2(I{\times}\Omega;\R^{d\times d})}^{}
\big\|\overlineStau\big\|_{L^2(I;H^1(\Omega;\R^{d\times d}))}^{}\le C,
    \\&\label{est:3}
    \big\|%\mtdEpk
    \overlineEptau\big\|_{L^\infty(I;H^1(\Omega;\R^{d\times d}))}^{}\le C\ \text{ and }\ \:
    %\mtdk\Epk\|_{L^2(I;H^1(\Omega;\R^{d\times d}))}^{}\le C,
       \Big\|\pdt\Eptau+(\overline\vv_\tau{\cdot}\nabla)\overlineEptau\Big\|_{L^2(I\times\Omega;\R^{d\times d})}^{}\le C,
    %\ \ \ \text{ and }\ \ \ \|\mtdEp\|_{L^2(I{\times}\Omega;\R^{d\times d})}^{}\le C,
    \\&\label{est:4}
    \|\overline\alpha_\tau\|_{L^\infty(I;H^1(\Omega;\R^\ell))}^{}\le C
   % \ \ \ \text{ and }\ \ \ \ 
,\ \ \ 
\Big\|%\mtd\alpha
    \pdt{\alpha_\tau}{+}(\overline\vv_\tau{\cdot}\nabla)\overline\alpha_\tau\Big\|_{L^2(I{\times}\Omega;\R^\ell)}^{}\!\le C, \ \text{ and }\ \
    \TTT\Big\|\pdt{\alpha_\tau}\Big\|_{L^2(I{\times}\Omega;\R^\ell)}^{}\!\le\frac{C}{\sqrt[4]{\tau}}\,,
\\&\label{est:4+}
    \|\overline\zz_\tau\|_{L^\infty(I;L^2(\Omega))}^{}\le C \ \ \ \ \ \text{ and }\ \ \ \ \ \|\overline\mu_\tau\|_{L^2(I;H^1(\Omega))}^{}\le C.
\end{align}
Actually, the estimate \eqref{est:1} is due to the Korn inequality.
%, while \eqref{est:2} has been obtained from
%%$\nabla\mtdEe=\nabla\EE(\vvk)-\nabla\mtdEp
%$\nabla\mtdEek=\nabla\EE(\vvk)-\mtdk\Epk$
%by using  %\eqref{eq:1b} together with \eqref{est:1} and \eqref{est:3}.
Moreover, from the calculus
\begin{align}\nonumber
\nabla{\bm S}=\nabla\partial_{\Ee}^{}\varphi(\Ee,\alpha,\zz)&=
\partial_{\Ee\Ee}^2\varphi(\Ee,\alpha,\zz)\nabla\Ee
%\\&\ \ \
+\partial_{\Ee\alpha}^2\varphi(\Ee,\alpha,\zz)\nabla\alpha
%\\&
+\partial_{\Ee\zz}^2\varphi(\Ee,\alpha,\zz)\nabla\zz
\ \ \text{ and}
%\,,\label{calculus-for-diffusion-Ee}
\\\nonumber
\nabla\mu=\nabla\partial_{\zz}^{}\varphi(\Ee,\alpha,\zz)&=
\partial_{\Ee\zz}^2\varphi(\Ee,\alpha,\zz)\nabla\Ee
%\\&\ \ \
+\partial_{\alpha\zz}^2\varphi(\Ee,\alpha,\zz)\nabla\alpha
%\\&
+\partial_{\zz\zz}^2\varphi(\Ee,\alpha,\zz)\nabla\zz\,,
%\label{calculus-for-diffusion-zz}
\end{align}
%\eqref{calculus-for-diffusion-Ee},
we can see
\begin{align}\nonumber
  \bigg(\!\!\begin{array}{c}\nabla\overlineEetau\\\nabla\overline\zz_\tau\end{array}\!\!\!\bigg)
  =\bigg(\!\!\begin{array}{cc}
  \partial_{\Ee\Ee}^2\varphi(\overlineEetau,\overline\alpha_\tau,\overline\zz_\tau)
&\partial_{\Ee\zz}^2\varphi(\overlineEetau,\overline\alpha_\tau,\overline\zz_\tau)\\\partial_{\Ee\zz}^2\varphi(\overlineEetau,\overline\alpha_\tau,\overline\zz_\tau)&\partial_{\zz\zz}^2\varphi(\overlineEetau,\overline\alpha_\tau,\overline\zz_\tau)
  \end{array}\!\!\bigg)^{-1}
\bigg(\!\!\begin{array}{c}\nabla\overlineStau
-
\partial_{\Ee\alpha}^2\varphi(\overlineEetau,\overline\alpha_\tau,\overline\zz_\tau)\nabla\overline\alpha_\tau\\
\nabla\overline\mu_\tau-
\partial_{\alpha\zz}^2\varphi(\overlineEetau,\overline\alpha_\tau,\overline\zz_\tau)\nabla\overline\alpha_\tau
\end{array}\!\!\!\bigg)
\,.
%\label{calculus-for-diffusion-Ee+}
\end{align}
From this, by using also the assumption \eqref{ass:1} which implies
boundedness of the inverse of the Hessian $\partial_{(\Ee,\zz),(\Ee,\zz)}^2\varphi$, we can still read the estimate
\begin{align}
\label{est:2++}
    %\|\Eek\|_{L^2(I;H^1(\Omega;\R^{d\times d}))}^{}\le C,
    \big\|\overlineEetau\big\|_{L^2(I;H^1(\Omega;\R^{d\times d}))
    }^{}\le C\ \ \ \ \text{ and}\ \ \ \ \ \big\|\overline\zz_\tau\big\|_{L^2(I;H^1(\Omega))
    }^{}\le C\,.
\end{align}\end{subequations}

%Furthermore, by the embedding $H^2(\Omega)\subset L^\infty(\Omega)$ for $d\le3$,
%we have $\vvk\in L^2(I;L^\infty(\Omega;\R^d))$; here we use that $\vvk(t)$
%ranges the space $H^2(\Omega;\R^d)$ so that the Bochner measurability
%of $\vvk$ into the nonreflexive space $L^\infty(\Omega;\R^d)$ indeed holds,
%although
%%while the Bochner measurability can be violated but
%it is not essential for
%the \SSS a priori \EEE estimation.
From the $L^\infty(I;L^2(\Omega))$-estimates of the gradients of
%$\Ee$,
$\overlineEptau$ and $\overline\alpha_\tau$
%and $\overline\zz_\tau$ already obtained in (\ref{est}c--e),
(\ref{est}c,d),
%\COMMENT{HERE I SEE A PROBLEM:\\
%  IT IS NOT CLEAR TO ME HOW $\nabla\bm{\varPi}_k
%=\nabla\pdt{}\Epk+\nabla((\vv\cdot\nabla)\Epk)
%\in L^2(L^2)$ CAN GIVE AN INFORMATION ABOUT $\nabla\Epk$???? THEN ALSO \eqref{est+} IS NOT CLEAR, AND THE STRONG CONVERGENCE OF $\Eek$ !!}
we can then estimate also
%%$\pdt{}\Eek=\EE(\vvk)-{\bm{\varPi_k}}
%%%-\mtdEek
%%-P_k((\bm v_k\cdot\nabla)\Eek)$, cf.\ \eqref{eq:1b+app},
%$\pdt{}\Eptau=\mtdk\Epk-P_k((\bm v_k\cdot\nabla)\Epk)$, 
%$\pdt{}\alpha\TTT_\tau^{k-1}\EEE=\mtd\alpha\TTT_\tau^{k-1}\EEE-((\bm v_k\cdot\nabla)\alpha\TTT_\tau^{k-1}\EEE)$, and
%$\pdt{}\zz\TTT_\tau^{k-1}\EEE)=\mtd\zz\TTT_\tau^{k-1}\EEE)-(\bm v_k\cdot\nabla\zz\TTT_\tau^{k-1}\EEE))=
%{\rm div}_k(\mathbb M(\zz\TTT_\tau^{k-1}\EEE))\nabla\mu_k)-P_k(\bm v_k\cdot\nabla\zz\TTT_\tau^{k-1}\EEE))$ \COMMENT{with ${\rm div}_k$ .......STILL TO DEFINE} as:
%, we obtain also\COMMENT{$d=3$ ....HERE TO BE CONTINUED}
\begin{subequations}\label{est+}\begin{align}
%\label{est:5}&\Big\|\pdt\Eek\Big\|_{L^2(I{\times}\Omega;\R^{d\times d})}^{}\le C,\\
&\label{est:6}\Big\|\pdt\Eptau\Big\|_{L^{4/3}(I{\times}\Omega;\R^{d\times d})}^{}\le C,%\COMMENT{\ MAYBE\ A\ LOWER\ EXPONENT}
\\&\label{est:7}\Big\|\pdt{\alpha_\tau}\Big\|_{L^{4/3}(I{\times}\Omega;\R^\ell)}^{}\le C,%\COMMENT{\ MAYBE\ A\ LOWER\ EXPONENT}
\\&\nonumber\Big\|\pdt{\zz_\tau}\Big\|_{L^2(I;H^1(\Omega)^*)}^{}
=\sup_{\|\widetilde\mu\|_{L^2(I;H^1(\Omega))}\le1}\int_0^T\!\!\!\bigg(\int_\Omega
\big(\mathbb M(\underline\alpha_\tau,\underline\zz_\tau)\nabla\overline\mu_\tau
-\overline\vv_\tau\big){\cdot}
\nabla\widetilde\mu
\\[-1.em]&\hspace{19em}-({\rm div}\overline\vv_\tau)\overline\zz_\tau\widetilde\mu\,\d x
+\int_\Gamma h\widetilde\mu
\,\d S\bigg)\d t
\label{est:8}\le C.
\intertext{For (\ref{est+}a,b), we used
  $\overline\vv_\tau\in L^\infty(I;L^2(\Omega;\R^d))\,\cap\,L^2(I;H^1(\Omega;\R^d))
  \subset L^4(I{\times}\Omega;\R^d)$ so that certainly
  $(\overline\vv_\tau{\cdot}\nabla)\overlineEptau\in L^{4/3}(I{\times}\Omega;\R^{d\times d})$ and 
  $(\overline\vv_\tau{\cdot}\nabla)\overline\alpha_\tau\in
  L^{4/3}(I{\times}\Omega;\R^\ell)$.
  Moreover, by $\pdt{}\Eetau=
\EE(\overline\vv_\tau)-\pdt{}\Eptau-(\overline\vv_\tau{\cdot}\nabla)\overlineEptau
  +k_{\rm e}\Delta\overlineStau
  -(\overline\vv_\tau{\cdot}\nabla)\overlineEetau$, cf.\ \eqref{eq:1b-disc+},
  we have also}
  &\label{est:7+}\Big\|\pdt\Eetau\Big\|_{L^2(I;H^1(\Omega;\R^{d\times d})^*)+L^{4/3}(I{\times}\Omega;\R^{d\times d})}^{}\le C.%\COMMENT{.... ONLY\ A\ SEMINORM!!}
\intertext{By comparison}
%$\mtd{\vvk}=(\operatorname{div}{\boldsymbol S}+{\bm f})/\varrho\in................$ so that, in view of $\pdt{}\vvk=\mtd{\vvk}-(\vvk\cdot\nabla)\vvk$
\nonumber&\pdt{\vv_\tau}=\frac{\operatorname{div}(\overlineStau
%\partial_{\Ee}\varphi(\overlineEetau,\overline\alpha_\tau,\overline\zz_\tau)
  +\overline{{\boldsymbol S}}_{\rm str,\tau}
+k_{\rm v}\EE(\overline\vv_\tau)
%  -{\rm div}(k_{\rm v}\nabla\EE(\overline\vv_\tau))
+\overline{\bm f}_\tau)}{\density}-(\overline\vv_\tau{\cdot}\nabla)\overline\vv_\tau
  -\frac12(\operatorname{div}\overline\vv_\tau)\overline\vv_\tau
  \intertext{ with $\overline{{\boldsymbol S}}_{\rm str,\tau}$ the piecewise constant
    interpolant of the structural stress, cf.\ \eqref{eq:1a-disc+} and
\eqref{structural-stress+},
  % and of $\operatorname{sym}\nabla{\vvk}=\mtdEe+\mtdEp$
%due to \eqref{eq:1b},
we have also}
  &\label{est:9}
\Big\|\pdt{\vv_\tau}\Big\|
%...............\le C\,,  \ \ \\Big\|\frac{\vvk{-}\vvkk}\tau\Big\|
_{L^2(I;H^3(\Omega;\R^d)^*)}\le C.%\COMMENT{.... ONLY\ A\ SEMINORM!!}
  %\ \ \ \text{ and }\ \ \ \|{\vvk}\|_{L^2(I;H^1(\Omega;\R^d))}^{}\le C.
\intertext{
%The latter estimate in \eqref{est:8} is due to the Korn inequality.
Here we used that, by \eqref{est:3},
$\nabla\overlineEptau\boxtimes\nabla\overlineEptau
-\frac12|\nabla\overlineEptau|^2{\boldsymbol I}
\in L^\infty(I;L^1(\Omega;\R^{d\times d}))$
and similarly, by \eqref{est:4},
also $\nabla\overline\alpha_\tau\TTT\boxtimes\EEE\nabla\overline\alpha_\tau
-\frac12|\nabla\overline\alpha_\tau|^2{\boldsymbol I}
\in L^\infty(I;L^1(\Omega;\R^{d\times d}))$, and also that $\varrho$
is assumed constant.
Also, for the limit passage in \eqref{VI-disc}, we need the estimates}
&\|\Delta\overlineEptau\|_{L^2(I{\times}\Omega;\R^{d\times d})}^{}\le C
  \qquad\text{and}\qquad
\|\Delta\overline\alpha_\tau\|_{L^2(I{\times}\Omega;\R^\ell)}^{}\le C\,,
\label{est-Delta}\end{align}\end{subequations}
which can be seen by comparison from (\ref{eq-disc+}c,d). 

\medskip\noindent{\it Step 4. (Convergence)}:
By the Banach selection principle, we obtain a subsequence
converging weakly* with respect to topologies indicated in
\eqref{est} and \eqref{est+}. Moreover, we now prove also the strong
convergence
\begin{subequations}\label{strong-nabla}\begin{align}
    &&&\nabla\overlineEptau\to\nabla\Ep&&\hspace*{-2em}\text{strongly in }\
  L^2(I{\times}\Omega;\R^{d\times d\times d})\,\text{ and}
\\&&&\nabla\overline\alpha_\tau\to\nabla\alpha&&\hspace*{-2em}\text{strongly in }\ L^2(I{\times}\Omega;\R^{d\times\ell})\,.
%\\&&&\nabla\overline\zz_\tau\,\to\,\nabla\zz&&\hspace*{-2em}\text{strongly in }\ L^2(I{\times}\Omega;\R^d)\,.
\end{align}\end{subequations}
%Here we use $\zeta$ bounded \COMMENT{???}
%and a growth restriction of %$\partial_{\Ep}\varphi$
%$\partial_\alpha^{}\varphi$ and of $\partial_\zz^{}\varphi$.
To prove it, we take sequences $\{\widetildeEptau\}_{\tau>0}^{}$ and
$\{\widetilde\alpha_\tau\}_{\tau>0}^{}$
%, and $\{\widetilde\zz_\tau\}_{\tau>0}^{}$
%valued in the respective finite-dimensional
%subspaces (so that $\Epk-\widetildeEpk$, $\alpha\TTT_\tau^{k-1}\EEE-\widetilde\alpha\TTT_\tau^{k-1}\EEE$, and
%$\zz\TTT_\tau^{k-1}\EEE)-\widetilde\zz\TTT_\tau^{k-1}\EEE)$ are legitimate test functions for the Galerkin
%approximation)
piecewise constant in time with respect to the partition with the time step
$\tau$ and, for $\tau\to0$, converging strongly
%(which is possible due to density of the union of these finite-dimensional
%subspaces in the strong topology)
towards $\Ep$ and $\alpha$,
%and $\zz$,
respectively. Then we can see that
\begin{subequations}\label{strong-nabla}\begin{align}\nonumber
& \int_0^T\!\!\!\int_\Omega k_{\rm p}|\nabla(\overlineEptau{-}\widetildeEptau)|^2\,\d x\d t
  =-\int_0^T\!\!\!\int_\Omega\bigg(\partial_{\mtd{\Ep}}\zeta\Big(\underline\alpha_\tau,\underline\zz_\tau;\pdt{\Eptau}+(\overline\vv_\tau{\cdot}\nabla)\overlineEptau,\pdt{\alpha_\tau}+(\overline\vv_\tau{\cdot}\nabla)\overline\alpha_\tau\Big)
  \\[-.3em]&
  %\nonumber
  \hspace{17em} +
  %\partial_{\Ee}^{}\varphi(\overlineEetau,\overline\alpha_\tau,\overline\zz_\tau)
  \overlineStau\bigg){:}(\overlineEptau{-}\widetildeEptau)
  %\\[-.3em]&\hspace{21em}
  +k_{\rm p}\nabla\widetildeEptau%\widetildeEpk
  \Vdots\nabla(\overlineEptau{-}\widetildeEptau)\,\d x\d t\to0%\,,
  \label{strong-nabla-Pi}
  \intertext{and similarly}
  %\\
  \nonumber
  &\int_0^T\!\!\!\int_\Omega k_{\rm a}|\nabla(\overline\alpha_\tau{-}\widetilde\alpha_\tau)|^2\,\d x\d t
  =-\int_0^T\!\!\!\int_\Omega
\bigg(\partial_{\mtd\alpha}\zeta\Big(\underline\alpha_\tau,\underline\zz_\tau;\pdt{\Eptau}+(\overline\vv_\tau{\cdot}\nabla)\overlineEptau,\pdt{\alpha_\tau}+(\overline\vv_\tau{\cdot}\nabla)\overline\alpha_\tau\Big)
\\[-.3em]&\hspace{7em}+\partial_\alpha^{}\varphi(\overlineEetau,\overline\alpha_\tau,\overline\zz_\tau)+\sqrt\tau\pdt{\alpha_\tau}\bigg)\cdot(\overline\alpha_\tau{-}\widetilde\alpha_\tau)
+k_{\rm a}\nabla\widetilde\alpha_\tau:\nabla(\overline\alpha_\tau{-}\widetilde\alpha_\tau)\,\d x\d t
  \to0\,.
\label{strong-nabla-alpha}
\end{align}\end{subequations}
Here we used \eqref{ass:1} so that $\partial_{\Ee}^{}\varphi(\overlineEetau,\overline\alpha_\tau,\overline\zz_\tau)$ and $\partial_\alpha^{}\varphi(\overlineEetau,\overline\alpha_\tau,\overline\zz_\tau)$ are bounded in
the respective $L^{6/5+\epsilon}(I{\times}\Omega)$-spaces while
$\overlineEptau-\widetildeEptau\to0$ and $\overline\alpha_\tau-\widetilde\alpha_\tau\to0$ strongly in $L^{6-\epsilon}(I{\times}\Omega;\R^{d\times d})$ and
$L^{6-\epsilon}(I{\times}\Omega;\R^\ell)$, respectively; this is due to the Aubin-Lions theorem, relying on \eqref{est:4+}
with \eqref{est:7}. In \eqref{strong-nabla-Pi}, we used
 that
$\partial_{\mtd\Ep}\zeta(\underline\alpha_\tau,\underline\zz_\tau;\pdt{}\Eptau+(\overline\vv_\tau{\cdot}\nabla)\overlineEptau,\pdt{}\alpha_\tau+(\overline\vv_\tau{\cdot}\nabla)\overline\alpha_\tau)$ is bounded
in $L^2(I{\times}\Omega;\R^{d\times d})$.
Similarly, in \eqref{strong-nabla-alpha}, we used that
$\partial_{\mtd\alpha}\zeta(\underline\alpha_\tau,\underline\zz_\tau;\pdt{}\Eptau+(\overline\vv_\tau{\cdot}\nabla)\overlineEptau,\pdt{}\alpha_\tau+(\overline\vv_\tau{\cdot}\nabla)\overline\alpha_\tau)$ is bounded
in $L^2(I{\times}\Omega;\R^\ell)$ and, moreover, that
$\|\sqrt\tau\pdt{}\alpha_\tau\|_{L^2(I\times\Omega;\R^\ell)}^{}
=\mathscr{O}\TTT(\EEE\sqrt[\TTT4\EEE]\tau)\to0$ due to \TTT the last estimate in
%\eqref{est:7}
\eqref{est:4}. \EEE

Based on the estimates \eqref{est:2++} and \eqref{est:7+}, we have
$\overlineEetau\to\Ee$ strongly in $L^2(I;L^{6-\epsilon}(\Omega;\R^{d\times d}))$
due to the Aubin-Lions theorem, generalized for functions whose
time-derivatives are measures as in \cite[Cor.7.9]{Roub13NPDE}. By the
interpolation with the estimate in $L^\infty(I;L^{2}(\Omega;\R^{d\times d}))$, we
have the strong convergence even in a smaller space, e.g.\ in
$L^4(I;L^{3-\epsilon}(\Omega;\R^{d\times d}))$.
Thanks to the growth condition \eqref{ass:1-} from which we have also
$|\varphi(\Ee,\alpha,\zz)|\le (1+|\Ee|^{5/2-\epsilon}\!
    +|\alpha|^{4-\epsilon}\!+|\zz|^{4-\epsilon})/\epsilon$, we can see that
    $\varphi(\overlineEetau,\overline\alpha_\tau,\overline\zz_\tau)$ converges
 strongly in $L^{6/5-\epsilon}(I{\times}\Omega)$.
 Taking into account also \eqref{strong-nabla}, we obtain the convergence in
 the structural stress \eqref{structural-stress+}, namely
 $\overline{{\boldsymbol S}}_{\rm str,\tau}\to {\boldsymbol S}_{\rm str}$
 strongly in $L^1(I{\times}\Omega;\R^{d\times d})$ with
${\boldsymbol S}_{\rm str}$ from \eqref{structural-stress}.
 Thus, noting that $\overlineStau=
 \partial_{\Ee}^{}\varphi(\overlineEetau,\overline\alpha_\tau,\overline\zz_\tau)$
 converges even strongly in $L^2(I{\times}\Omega;\R^{d\times d})$ due to the
 growth condition \eqref{ass:1-}, we can pass to the limit in the momentum
 equation \eqref{eq:1a-disc+}. The limit passage in \eqref{eq:1b-disc+} is
 similar.

By the proved strong convergence of
$\overlineEetau\to\Ee$, we can pass to the limit in the nonlinear terms
$\partial_\Ee^{}\varphi(\overlineEetau,\overline\alpha_\tau,\overline\zz_\tau)$
and
$\partial_\alpha^{}\varphi(\overlineEetau,\overline\alpha_\tau,\overline\zz_\tau)$
flow rule, i.e.\ in the variational inequality %\eqref{eq:4-disc+}
\eqref{VI-disc}, and in the terms
$\partial_{\zz}^{}\varphi(\overlineEetau,\overline\alpha_\tau,\overline\zz_\tau)$
and $\mathbb M(\underline\alpha_\tau,\underline\zz_\tau)$
in the
%Cahn-Hilliard \COMMENT{NOW NOT 4th-order EQS.!! MAYBE MAYBE JUST diffusion equation???}
diffusion equation \eqref{eq:5-disc+}, too.

Let us also note that the convexifying term in \eqref{eq:4-disc+} vanishes
in the limit due to the estimate \eqref{est:7} because obviously
$\|\sqrt\tau\pdt{}\alpha_\tau\|_{L^2(I\times\Omega;\R^\ell)}^{}
=\mathscr{O}\TTT(\EEE\sqrt[\TTT4\EEE]\tau)\to0$, as used already before in
\eqref{strong-nabla-alpha}.

Eventually, from \eqref{est-Delta} and \eqref{est:2++}, we also %give
\TTT obtain \EEE the $L^2(I{\times}\Omega)$-information
about $\Delta\bm\varPi$, $\Delta\alpha$, and $\nabla\zz$.
%, which is a part of Definition~\ref{def}.
%\end{proof}

$\hfill\Box$

%\bigskip

%\COMMENT{ONE (NOT MUCH IMPRESSIVE) OPTION}
%put the spherical hyperstress ${\rm div}(k\nabla({\rm div}\vv){\bm I})$ into \eqref{eq:1a} -- assuming $\Ep$ isochoric, so that ${\rm div}\vv=\mtdEe$ and controlling the total spherical strain is acceptable.

%
%so that $\nabla({\rm div}\vv)\in L^2(I\times\Omega;\R^d)$ so that $({\rm div}\vv)$ will converge strongly by the Aubin-Lions theorem

%$\varphi$ and $\zeta$ quadratic in $\Ee$ and $\mtd{\Ee}$ so that the weak limit
%passage in \eqref{eq:1a} will work

%\COMMENT{... limit passage via weak convergence and Minty's trick
%  if $\varphi(\cdot,\alpha)$ is convex seems difficult since the test
%  by $\Ee$ is needed while only the tests by $\pdt\Ee$ or $\mtd{\Ee}$ work well}

%\COMMENT{... IT SEEMS THAT anyhow the viscous hyperstresses are needed
%  in the form $\Delta\mtd{\Ee}$, $\Delta\mtd{\Ep}$, and $\Delta\mtd\alpha$}

%\begin{remark}
%  Instead of testing \eqref{eq:1a} by $\vv_k$ (in the Faedo-Galerkin
%  approximation) and using the calculus \eqref{test-of-convective},
%  we can use the form \eqref{eq:1a+} and the calculations from
%  Remark~\ref{rem-Bernoulli}.
%\end{remark}

\section{Concluding remarks}
%        ~~~~~~~~~~~~~~~~~~

We close this paper with several remarks, outlining some concrete
examples, expansions, or comments to the used analysis.

\begin{remark}[Example for a semi-covex $\varphi$.]\label{rem-example-energy}
The so-called (weakened) semi-convexity \eqref{semi-convex}
%{ass:1-}
is not in conflict with usual damage models and, when combined with
{\it Biot's poroelasticity}, it allows for models like
\begin{align}\label{varphi}
  \varphi(\Ee,\alpha,\zz)=\frac{dK}2|{\rm sph}\Ee|^2
+\frac M2|\beta{\rm tr}\,\Ee{-}\zz+\zz_{\rm eq}|^2
+G(\alpha)\frac{|{\rm dev}\Ee|^2}{1{+}\epsilon|{\rm dev}\Ee|^2}
+G_0|{\rm dev}\Ee|^2+\phi(\alpha)\,,\end{align}
where  $\zz_{\rm eq}$ is a given equilibrium concentration, ``sph'' denotes
the spherical part (recall that
${\rm sph}\Ee=\Ee{-}{\rm dev}\Ee=({\rm tr}\Ee)\bm I/d$) with $K$ the bulk
modulus, ``tr'' denotes the trace, and ``dev'' the deviatoric part
with the shear modulus $G:\R^\ell\to\R$ non-negative smooth satisfying $G_i'(...,0,...)=0=G_i'(...,1,...)$ for $i=1,...,\ell$, which ensures
that each $\alpha_i$ takes values in the interval $[0,1]$ as usually
requires in damage/breakage type models. Further
parameters $K$, $M$, and $\beta$ in \eqref{varphi} have the meaning
of the bulk modulus, Biot's modulus, and Biot coefficient, respectively, while
$G_0>0$ is just small regularizing modulus not subjected to damage and
ensuring coercivity. This is the classical Biot model for a saturated
fluid flow in poroelastic media \cite{Biot41GTTC}.
Note that the second derivatives of the $G(\alpha)$-term
are bounded so that \eqref{ass:1-} holds.
For a convexification by a quadratic form in $(\Ee,\zz)$ see
\cite{Roub19CTDD} which deals with a non-convective variant and
which would be here more difficult.
Actually, the Biot ansatz \eqref{varphi} gives the chemical potential
$\mu= M(\beta{\rm tr}\Ee{-}\zz)$, meaning a pressure and then the flux in
the Fick diffusion turns rather to the {\it Darcy law}.
%\COMMENT{STILL ... Darcy diffusion.....}
The last term in \eqref{varphi} creates a driving force for healing of
damage. Together with the $\Delta\alpha$ in the damage flow-rule \eqref{eq:4},
it enables to model the Ambrosio-Tortorelli-type {\it phase-field
fracture}; actually, the standard choice is $G(\alpha)=G_1\alpha^2$,
$G_0=k_{\rm a}^2$ and $\phi(\alpha)=(1{-}\alpha)^2/(2k_{\rm a})$ with 
$k_{\rm a}>0$ from \eqref{eq:4} assumed small.
\COMMENT{REFERENCES?}
\end{remark}

\begin{remark}[Energy conservation.]\COMMENT{NEW REMARK:}
The energy balance \eqref{energy} is only formal and its rigorous proof
needs to legitimate the test used in \eqref{eq:12}--\eqref{eq:18}.
This does not seem easily possible, however, and a regularization
of the model seems necessary. More specifically, a higher-order viscosity
of the type
${\rm div}\FFF^2\EEE(k_{\rm v}\FFF|\nabla\EE(\vv)|^{p-2}\EEE\nabla\EE(\vv))$
\FFF for $p>d$ \EEE together with the viscous
%Cahn-Hilliard
\FFF variant of the \EEE diffusion
$\mu=\partial_\zz\varphi(\Ee,\alpha,\zz)+\epsilon\mtd{}\zz%-\epsilon\Delta\zz
$ \FFF with some (presumably small) modulus $\epsilon>0$ (with the physical
dimension Pa\,s=J\,s/m$^3$) \EEE
would help\FFF, cf.\ \cite{Roub??SPTC} or \cite[Sect.8]{Roub??TCMP}. \EEE
This would \FFF only \EEE
%enhance in particular the structural stress and
make
%a lot of calculations and
\FFF some \EEE arguments \FFF a bit \EEE more complicated
%. On the other hand, it would
\FFF and \EEE open \FFF a possibility for \EEE an expansion of the model
towards full thermodynamics by %closing
\TTT completing \EEE
it by the heat-transfer equation. As for the analysis, first the limit
passage in the mechanical part using also the strong convergence
  of temperature by the Aubin-Lions compactness theorem should be done,
  %and then,
  \FFF followed by the \EEE strong convergence of the dissipation rate,
  %by using that $\pdt\vv$ is in duality with $\vv$ \COMMENT{TO CHECK: the strucural stress in $L^\infty(I;L^3(\Omega;\R^{d\times k}))$  because  $|\nabla\alpha|^2\in L^\infty(I;L^3(\Omega;\R^{d\times k}))$ and $|\nabla\zz|^2\in L^\infty(I;L^3(\Omega;\R^d))$    by $H^2$-regularity}
  and % eventually
  \FFF finished by the \EEE convergence in the heat-transfer equation. We refer
  to \cite[Chap.\,8]{KruRou19MMCM} \TTT or also e.g.\ \cite{Roub??SPTC} \EEE
for the technical details.
\end{remark}

\begin{remark}[Staggered time discretisation.]
  One could think about a fractional-step splitting
  (also known as a staggered) time discretisation
  to  decouple $(\vvk,\Eek,\Epk)$ from $\alphak$ and from
  $(\zzk,\muk)$ in order to allow for a separately convex
  $\varphi$.
  %with a more general growth, when consider the less general
  %dissipation potential  $\zeta(\alpha,\zz;\mtd\Ep,\mtd\alpha)=
  %\zeta_{\rm p}(\alpha,\zz;\mtd\Ep)+\zeta_{\rm a}(\alpha,\zz;\mtd\alpha)$.
  This usually works efficiently,
although here it would lead to a coupled scheme through the structural
stress but, more important, here there would be troubles with modification
of the calculus \eqref{calc-energy-pressure}. This is the reason that we used
the fully implicit time discretisation \eqref{eq-disc}.
\end{remark}

  \begin{remark}[Galerkin method.]
 In our convective model, the Galerkin approximation (i.e.\ the space discretisation instead of the time discretisation \eqref{eq-disc}) would face serious
technical difficulties because testing by convective time derivatives
which do not comply with finite-dimensional spaces used for the Galerkin
method and the sophisticated calculus like %\eqref{calc-to-Korteweg}
\eqref{calc-plastic}--\eqref{calc-alpha} or
\eqref{calc-energy-pressure} would not be legal.
%Also, the $H^2$-regularity can be used only
%after making the limit passage towards the continuous problem and
%the strong convergence of $\nabla\Ep$, $\nabla\alpha$, and $\nabla\zz$
%needed for the nonlinear structural stress would have to be proved
%directly. For the former difficulty
Therefore, the implementation of this, usually
very efficient technique seems problematic here. 
%and the estimates \eqref{est:10} can be used 
\end{remark}

%\begin{remark}\COMMENT{TO MODIFY}
%When using a time discretisation instead of the space discretisation,
%the $H^2$-regularity and the estimates \eqref{est:10} can be used already
%on the discrete level. The strong convergence of $\nabla\alpha$'s and
%$\nabla\zz$'s would then follow simply by the Aubin-Lions compact
%embedding theorem, omitting thus the calculations \eqref{strong-nabla-alpha}
%and \eqref{strong-nabla-zz} but the assumption \eqref{ass:1} is to be
%strengthened for $|\varphi'(\Ee,\alpha,\zz)|\le C(1+|\Ee|^{5}+|\alpha|^3+|\zz|^3)$  so that $\varphi_{(\alpha,\zz)}'(\Ee,\alpha,\zz)$ would be controlled in
%$L^2(I{\times}\Omega;\R^{k+1})$ and $\Omega$ should be smooth or convex.
%\end{remark}

%\begin{remark}\COMMENT{TO MODIFY}
%  In geophysical simulations ................, usually a time-evolving shift
%  of the boundaries of the computational domain $\Omega$ is prescribed,
%  which would correspond to prescribe a nonvanishing $\vvn$ in \eqref{BC1}.
%  This, unfortunately, would bring essential difficulties in our
%  analysis. More specifically, the estimation of the last term in
%  \eqref{test-of-convective+} might be delicate. In principle, one may
%  consider a condition like $\vvn=\vv_{\rm prescribed}/(1+\epsilon|\vv|^2)$,
%  which would made this term apriori bounded, but then other difficulties
%  would occur in ............................................
%\end{remark}

\begin{remark}[More general stored energies.]
The stored energy $\varphi$ is often %not quadratic (or even possibly
considered not convex in geophysical applications, as devised in
\cite{LyaMya84BECS} and used e.g.\ in
\cite{LyaBeZ14CDBF,LyaHam07DWFF,LyHaBZ11NLVE}.
  This brings, however, technical difficulties in analysis. In particular, 
  it violates the assumption \eqref{ass:1} which is needed to control
  $\nabla\Ee$ which was used to obtain strong convergence in $\Ee$.
  And this strong is needed to pass to the limit in the nonlinear
  terms $\varphi(\cdot,\alpha,\zz)\bm I$ and
  %Passing through a monotone
  $\partial_{\Ee}^{}\varphi(\cdot,\alpha,\zz)$ in particular in such a
  nonconvex situation.
  % by the weak convergence should be based on testing
  %$\varphi_{\Ee}'(\Eek,\alpha\TTT_\tau^{k-1}\EEE,\overline\zz_\tau)$ by $\Eek$
  %...... In case of a nonmonotone $\varphi_{\Ee}'(\cdot,\alpha,\zz)$, ......
\end{remark}

 % \begin{remark}[Thermodynamical extension]\label{rem-thermo}
 %   The mechanical-energy conservation allows for possibilities of
 %   completion of the system by the heat generation and transfer and
 %   for making the dissipative mechanisms also temperature dependent.
%%    \COMMENT{HOW THERMODYNAMICS WOULD WORK NOW??}
 %   As for the analysis, first the limit passage in the mechanical part, using
 %   also the strong convergrence
 % of temperature by the Aubin-Lions compactness theorem should be done, and
 % then, strong convergence of the dissipation rate.
 % %by using that $\pdt\vv$ is in duality with $\vv$ \COMMENT{TO CHECK: the strucural stress in $L^\infty(I;L^3(\Omega;\R^{d\times k}))$  because  $|\nabla\alpha|^2\in L^\infty(I;L^3(\Omega;\R^{d\times k}))$ and $|\nabla\zz|^2\in L^\infty(I;L^3(\Omega;\R^d))$    by $H^2$-regularity}
 % and  eventually convergence in the heat-transfer equation. We refer
 % to \cite[Chap.\,8]{KruRou19MMCM} for the technical details.
%\end{remark}

\begin{remark}[Other phenomena involved.]
The Eulerian description opens a way for enhancement of the model by other phenomena which ultimately needs formulation in Eulerian configuration. In particular, it concerns gravity and magnetic fields. Also, a coupling with fluidic regions (in particular with the outer core of the Earth) is thus well facilitated. 
\end{remark}

\section*{Acknowledgements}
This research has been partially \TTT supported \EEE from the
CSF (Czech Science Foundation) project 19-04956S, the M\v SMT \v CR
(Ministry of Education of the Czech Rep.) project
CZ.02.1.01/0.0/0.0/15-003/0000493, and the institutional support RVO: 61388998 (\v CR). This research has also been partially supported by the Italian INdAM-GNFM (Istituto Nazionale di Alta Matematica-Gruppo Nazionale per la Fisica Matematica), the Grant of Excellence Departments, MIUR-Italy (Art.1, commi 314-337, Legge
232/2016), and the Grant ``Mathematics of active materials: from mechanobiology to smart devices'' (PRIN 2017, prot. 2017KL4EF3) funded by the Italian MIUR.

%\TTT
%\begin{acknowledgements}
%The authors are deeply thankful to two anonymous referees for very careful
%reading of the original version and many suggested improvements and
%mistake corrections. 
%\end{acknowledgements}
%\EEE

 %Authors must disclose all relationships or interests that 
 %could have direct or potential influence or impart bias on 
 %the work: 

% \section*{Conflict of interest}The authors declare that they have no conflict of interest.

%\bibliographystyle{plain}
%\bibliography{tr-gt-geo}

\end{document}